\documentclass{gtpart}   
\agtart

\usepackage{amscd,amsfonts,amssymb,amscd,amsmath,latexsym}

\usepackage[all]{xy}

\makeatletter
\renewcommand{\subsubsection}{\@startsection
{subsubsection}
{1}
{0mm}
{0mm}
{0mm}
{\normalfont\normalsize}}
\makeatother

\title{Landweber exact formal group laws and smooth cohomology theories}

\author{Ulrich Bunke} 
\givenname{Ulrich}
\surname{Bunke}
\address{NWF I - Mathematik \\Universit\"at Regensburg\\93040 Regensburg \\ Germany}
\email{ulrich.bunke@mathematik.uni-regensburg.de}
\urladdr{http://www.mathematik.uni-regensburg.de/Bunke/}

\author{Thomas Schick} 
\givenname{Thomas}
\surname{Schick}
\address{Mathematisches Institut\\Georg-August-Universit\"at
  G\"ottingen\\Bunsenstr.~3\\37073 G\"ottingen\\German} 
\email{schick@uni-math.gwdg.de}
\urladdr{http://www.uni-math.gwdg.de/schick}

\author{Ingo Schr\"oder} 
\givenname{Ingo}
\surname{Schr\"oder}
\address{Mathematisches Institut\\Georg-August-Universit\"at G\"ottingen\\Bunsenstr.~3\\37073 G\"ottingen\\Germany}
\email{ischroed@uni-math.gwdg.de}

\author{Moritz Wiethaup} 
\givenname{Moritz}
\surname{Wiethaup}
\address{Mathematisches Institut\\Georg-August-Universit\"at
  G\"ottingen\\Bunsenstr.~3\\37073 G\"ottingen\\Germany} 
\email{wiethaup@uni-math.gwdg.de}


\numberwithin{equation}{section}

 \newtheorem{theorem}{Theorem}[section] 
\newtheorem{prop}[theorem]{Proposition}
\newtheorem{lem}[theorem]{Lemma}
\newtheorem{ddd}[theorem]{Definition}
\newtheorem{kor}[theorem]{Corollary}

\newcommand{\hMU}{\hat{MU}}

\newcommand{\MM}{\mathbb{M}}
\newcommand{\bMG}{{\mathbf{MG}}}

\renewcommand{\lim}{{\tt lim\:}}
\newcommand{\colim}{{\tt colim}}
\renewcommand{\P}{{\mathbb{P}}}

\renewcommand{\det}{{\tt det}}

\newcommand{\Tor}{{\tt Tor}}

\newcommand{\End}{{\tt End}}

\newcommand{\im}{{\tt im}}

\newcommand{\clo}{{\tt clo}}

\newcommand{\Gr}{{\tt Gr}}

\newcommand{\id}{{\tt id}}

\newcommand{\nat}{\mathbb{N}}

\def\imath{{i}}

\newcommand{\pr}{{\tt pr}}

\begin{document}
\def\thedoi{10.2140/agt.2009.9.1751}
\def\thestartpage{1751}
\def\theendpage{1790}
\def\thepublicationyear{2009}
\def\thevolumenumber{9}
\received{24 September 2008}
\revised{15 July 2009}
\accepted{19 July 2009}
\def\publisheddate{26 September 2009}

\begin{abstract}
  The main aim of this paper is the construction of a smooth (sometimes called
  differential) extension $\widehat{MU}$ of the cohomology theory complex
  cobordism $MU$, using 
  cycles for $\widehat{MU}(M)$ which are essentially proper maps $W\to M$ with
  a fixed 
  $U $-structure and $U $-connection on the (stable) normal bundle of $W\to M$. 

  Crucial is that this model allows the construction of a product structure
  and of pushdown maps for this smooth extension of $MU$, which have all the
  expected properties.

  Moreover, we show, using the Landweber exact functor principle, that
  $\hat{R}(M):=\widehat{MU}(M)\otimes_{MU^*}R$ defines a multiplicative smooth
  extension of $R(M):=MU(M)\otimes_{MU^*}R$ whenever $R$ is a Landweber exact
  $MU^*$-module. An example for this construction is a new way to define
 a multiplicative smooth K-theory.
\end{abstract}

\maketitle

\section{Introduction}

Smooth (also called differentiable) extensions of generalized cohomology theories recently became an
intensively studied mathematical topic with many applications ranging from arithmetic geometry to string theory. Foundational contributions are \cite{MR827262}, \cite{MR1197353} (in the case of ordinary cohomology) and
\cite{MR2192936}. The latter paper gives among many other results  a general
construction of  smooth extensions in homotopy theoretic terms. For cohomology
theories based on geometric or analytic cycles there are often alternative
models. This applies in particular to ordinary cohomology whose smooth
extension has various different realizations (\cite{MR827262},
\cite{MR1423029}, \cite{MR1197353}, \cite{MR2179587}, \cite{MR2192936},  \cite{bunke-kreck}). The
  papers \cite{math.AT/0701077} or
\cite{bs2009} show that all
these realizations are isomorphic. 

 An example of a cycle model of a smooth extension of a generalized cohomology theory is the model of smooth  $K$-theory introduced  
  in \cite{bunke-2007}, see also \cite{hep-th/0011220}, \cite{MR1769477}. The present paper contributes geometric models of smooth extensions of cobordism theories, where the case of complex cobordism theory $MU$ is of particular importance. {
In \cite{bs2009} we obtain general results about uniqueness of smooth
extensions which in particular apply to smooth $K$-theory and   smooth complex
cobordism theory $\hMU$. In detail, any two smooth extensions of complex cobordism
theory or complex $K$-theory  which admit an integration along $\int\colon S^1\times M\to M$  
are isomorphic by a unique isomorphism compatible with $\int$. In case of multiplicative extensions the isomorphism is automatically multiplicative. Note that
the extension $\hMU$ constructed in the present paper has an integration and is multiplicative. }

{
We expect that our model $\hMU$ of the smooth extension of $MU$ is uniquely isomorphic to
the one given by \cite{MR2192936}. So far this fact can not immediately be deduced from the above uniqueness result since for the latter model the functorial properties of the integration map have not been developed yet in sufficient detail.
However, for the uniqueness of the even part we do not need the integration. Therefore in even degrees our extension $\hMU$ is uniquely isomorphic to the model in \cite{MR2192936}}.

An advantage of geometric or analytic models is that they allow the introduction
of additional structures like products, smooth orientations and integration maps with good properties.
These additional properties are fundamental for applications. In \cite[4.10]{MR2192936} methods for integrating
  smooth cohomology classes
were discussed, but further work will be required in order to turn these
ideas into constructions with good functorial properties.

In the case of smooth ordinary cohomology the product and the integration have been considered  in various places 
(see e.g. \cite{MR827262}, \cite{MR2179587}, \cite{MR1197353}) (here smooth orientations are ordinary orientations).
The case of smooth $K$-theory, discussed in detail in \cite{bunke-2007}, shows that
in particular the theory of orientations and integration is considerably more complicated for generalized
cohomology theories.

 In the present paper we construct a multiplicative extension of the {complex cobordism cohomology theory} $MU$. Furthermore, we introduce the notion of a smooth $MU$-orientation and develop the corresponding theory of integration.
The same ideas could be applied with minor modifications to other cobordism theories.

For an $MU^*$-module $R$ one can try to define a new cohomology theory
$R^*(X):=MU^*(X)\otimes_{MU^*}R$ for finite $CW$-complexes $X$. By
 Landweber's famous result \cite{MR0423332} this construction
works and gives a multiplicative complex oriented cohomology theory provided
$R$ is a ring over $MU^*$ which is in addition  \textit{Landweber exact}. In
Theorem \ref{hjgwzuwd} we observe that by the same idea 
one can define a multiplicative smooth extension $\hat R(X):=\hMU(X)\otimes_{MU^*}R$ of $R$.
It immediately follows that this smooth extension admits an integration for smoothly $MU$-oriented proper submersions.

In this way we considerably enlarge the class of examples of generalized cohomology theories which admit multiplicative extensions and integration maps. The construction can e.g. be applied to
Landweber exact elliptic cohomology theories \cite{MR1320998},
\cite{MR1235295} and complex $K$-theory \footnote{It is an interesting problem
  to understand explicitly  the relation with the model
  \cite{bunke-2007}. Note that we abstractly know that the
    smooth extensions are isomorphic by \cite{bs2009}.}.

In Section \ref{dfgs}   we review the main result of Landweber \cite{MR0423332} and 
the definition of a smooth extension of a generalized cohomology theory. We state the main result asserting 
the existence of a multiplicative smooth extension of $MU$ with orientations and integration.
Then we realize the  idea sketched above and construct a multiplicative smooth extension for
every Landweber exact formal group law.

In Section \ref{gzdqw} we review the standard constructions of cobordism theories
using homotopy theory on the one hand,  and cycles on the other. Furthermore, we review the notion of a genus.

In Section  \ref{wuzwqud} we construct our model of the multiplicative smooth extension of complex cobordism. Furthermore, we introduce the notion of a smooth $MU$-orientation and construct the integration map.

\bigskip

{\small Thomas Schick was {partially funded by the Courant Research Center ``Higher
  order  structures in Mathematics'' within the German initiave of
  excellence}. Ingo Schr\"oder and Moritz Wiethaup were {partially
    funded by DFG GK 535 ``Groups and Geometry''}.}

\section{The Landweber construction and smooth extensions}\label{dfgs}

\subsection{The Landweber construction }
\subsubsection{}

Let $X\mapsto MU^*(X)$ denote the multiplicative cohomology theory (defined on the category of $CW$-complexes) called complex cobordism.
We fix an isomorphism $MU^*(\C\P^\infty)\cong MU^*[[x]]$. The K\"unneth formula then gives
$MU^*(\C\P^\infty\times \C\P^\infty)\cong MU^*[[x,y]]$.

The tensor product of line bundles induces an $H$-space structure $\mu\colon \C\P^\infty\times \C\P^\infty\to  \C\P^\infty$.
Under the above identifications the map
$\mu^*\colon MU^*[[z]]\to MU^*[[x,y]]$ is determined by the element
$f(x,y):=\mu^*(z)\in MU^*[[x,y]]$. 

By a result of Quillen \cite{MR0253350} the pair $(MU^*,f)$ is a universal formal group law. 
This means that, given a commutative ring $R$ and a formal
group law $g\in R[[x,y]]$, there exists a unique
ring homomorphism $\theta\colon MU^*\to R$ such that $\theta_*(f)=g$.

\subsubsection{}

Let $R$ be a commutative ring over $MU^*$. Then one can ask if the functor
$X\mapsto MU^*(X)\otimes_{MU^*} R$ is  a cohomology theory on the category of finite $CW$-complexes. 
The result of Landweber \cite{MR0423332} determines necessary and sufficient conditions.
A ring which satisfies these conditions is called Landweber exact.

\subsubsection{}

Actually, Landweber shows a stronger result which is crucial for the present paper.
For any space or spectrum $X$ the homology $MU_*(X)$ has {the} structure of a comodule over the coalgebra
$MU_*MU$ in $MU^*$-modules. By duality, if $X$ is finite, then $MU^*(X)\cong MU_*(S(X))$ also has a comodule structure, where
$S(X)$ denotes the Alexander-Spanier dual (see \cite{MR0402720}) of $X$.

\begin{theorem}[Landweber \cite{MR0423332} ]\label{landweber}
 Let $M$ be a  finitely presented $MU^*$-module which has the structure of a comodule over $MU_*MU$, 
and consider a Landweber exact formal group law
$(R,g)$ so that in particular $R$ is a ring over $MU^*$. 
Then  for all $i\ge 1$ we have $\Tor^{MU^*}_i(M,R)=0$.
\end{theorem}

\subsection{Smooth cohomology theories}

\subsubsection{}\label{uisadusd}
{In the present subsection $B$ denotes a compact smooth manifold.}
  Let $N$ be a $\Z$-graded vector space over $\R$. We consider  a generalized cohomology theory $h$ with a  natural transformation of cohomology theories
$c\colon h(B)\to H(B,N)$,  where $H(B,N)$ is ordinary cohomology with
  coefficients in $N$. The natural universal example is given by $N:=h^*\otimes \R$, where $c$ is the canonical transformation. Let $\Omega(B,N):=\Omega(B)\otimes_\R N$, where $\Omega(B)$ denotes the smooth real differential forms on $B$. {Note that this definition only coincides with the corresponding definition of $N$-valued forms in \cite{bs2009} if $N$ is degree-wise finite-dimensional.}
By $dR\colon \Omega_{d=0}(B,N)\to H(B,N)$ we denote the de Rham map which associates to a closed form the corresponding cohomology class.
To a pair $(h,c)$ we associate the notion of a smooth extension $\hat h$. 
Note that manifolds in the present paper may have boundaries. 
\begin{ddd}\label{ddd556}
A smooth extension of the pair $(h,c)$ is a functor $B\to\hat h(B)$ from the category of compact smooth manifolds
to $\Z$-graded groups together with natural transformations
\begin{enumerate}
\item $R\colon  \hat h(B)\to \Omega_{d=0}(B,N)$ (curvature)
\item $I\colon  \hat h(B)\to  h(B)$ (forget smooth data)
\item $a\colon \Omega(B,N)/\im(d)\to  \hat h(B)$ (action of forms)\ .
\end{enumerate}
These transformations are required to satisfy the following axioms:
 \begin{enumerate}
\item The following diagram commutes
$$\xymatrix{\hat h(B)\ar[d]^R\ar[r]^I& h(B)\ar[d]^{c}\\
\Omega_{d=0}(B,N)\ar[r]^{dR}&H(B,N)}\ .$$
\item \begin{equation}\label{drgl}
R\circ a=d\ .
\end{equation}
\item $a$ is of degree $1$.
\item The sequence
\begin{equation}\label{exax}
h(B)\stackrel{c}{\to} \Omega(B,N)/\im(d)\stackrel{a}{\to} \hat h(B)\stackrel{I}{\to} h(B)\to 0\ .
\end{equation}
is exact.
\end{enumerate}
\end{ddd}

\subsubsection{}

If $h$ is a multiplicative cohomology theory, then one can consider a
$\Z$-graded ring $R$ over $\R$ and a  multiplicative transformation $c\colon h(B)\to
H(B,R)$. In this case we also talk about a multiplicative smooth extension $\hat h$ of  $(h,c)$.
\begin{ddd}\label{multdef1}
A smooth extension $\hat h$ of $(h,c)$ is called multiplicative, if $\hat h$ together with the transformations $R,I,a$ is a smooth extension of $(h,c)$, and in addition
\begin{enumerate}
\item $\hat h$ is a functor to $\Z$-graded rings,
\item $R$ and $I$ are multiplicative, 
\item $a(\omega)\cup x=a(\omega\wedge R(x))$ for $x\in \hat h(B)$ and $\omega\in\Omega(B,R)/\im(d)$. 
\end{enumerate}
\end{ddd}
 
\subsubsection{}

The first goal of the present paper is the construction of a multiplicative smooth extension of the pair
$(MU,c)$, where $c\colon MU^*(B)\to MU^*(B)\otimes_\Z\R\cong  H^*(B,MU\R)$ is the canonical natural transformation (see \ref{univcase}). The following theorem  is a special case of Theorem \ref{fundprop} which gives a construction of multiplicative smooth extensions of more general pairs $(MU,h)$.
\begin{theorem}\label{ex}
The pair $(MU,c)$ admits a multiplicative smooth extension.
\end{theorem}
The existence of a smooth extension also follows from \cite{MR2192936},
but there, no ring structure is constructed.

\subsubsection{}

In the present paper we consider smooth extensions of
  generalized cohomology  theories defined on the category of compact manifolds. The reason lies in the fact that we want to apply the Landweber exact functor theorem. If $R$ is a generalized complex oriented cohomology theory satisfying the wedge axiom to which the Landweber exact functor theorem applies, then for finite $CW$-complexes $X$
$$R^*(X)\cong MU^*(X)\otimes_{MU^*} R\ . $$ In general this
equality does not extend to infinite $CW$-complexes since the tensor product
on the right-hand side  does  not necessarily commute with infinite products. 

If one omits the compactness condition in the Definitions \ref{ddd556} and \ref{multdef1}, then one obtains the axioms for smooth  and multiplicative smooth extensions
defined on the category of all manifolds. If the  coefficients groups $R$ is degree-wise
finitely generated (see the corresponding remark in \ref{uisadusd}), then we obtain the same notion as in \cite{bs2009} 

Our construction of the smooth extension of the complex cobordism theory does not depend on any compactness assumption so that there is also a corresponding version of  Theorem \ref{ex} furnishing a multiplicative smooth extension of $(MU,c)$ defined on the category of all smooth manifolds.

\subsubsection{}\label{dhj}

We also introduce the notion of a smooth $MU$-orientation (Definition \ref{smoop}) of a proper submersion  $p\colon W\to B$ and define a push-forward
$  p_! \colon \hMU(W)\to \hMU(B)$ which refines the integration map
$p_!\colon MU(W)\to MU(B)$ (Definition \ref{cpush}).
In Subsection \ref{ghgtg} we show that
integration is compatible with the structure maps $a,R,I$ of the smooth extension, functorial,
compatible with pull-back and the product. We refer to this subsection and Theorem \ref{multi}
for further details. Integration maps play a fundamental role in applications
of generalized cohomology theories. This is the case e.g.~in the context of
T-duality, where we hope to eventually generalize our investigatons
\cite{MR2130624} to a setting in smooth cohomology.

\subsection{Smooth extensions for Landweber exact formal group laws}

\subsubsection{}

If $(R,g)$ is a Landweber exact formal group law, then we let
$R^*(X):=MU^*(X)\otimes_{MU^*} R$ denote the associated cohomology theory on finite $CW$-complexes.
We consider the pair $(R,c_R)$, where $c_R\colon R\to R\otimes_\Z \R=:R\R$ is the canonical map.

\begin{theorem}\label{hjgwzuwd}
If $(R,g)$ is a Landweber exact formal group law, then  
$(R,c_R)$ has a multiplicative smooth extension $\hat R$, given by $\hat{R}(B)=\hMU(B)\otimes_{MU^*}R$.
 \end{theorem}
\begin{proof}
We start with Theorem \ref{ex} which states that $(MU,c)$ has a multiplicative smooth extension. 
Since $\Omega^k(*)=0$ for $k\not=0$, $\Omega^0(*)\cong \R$, and $MU^{odd}=0$, the natural map
$\hMU^{ev}(*)\to MU^{ev}(*)$ is an  isomorphism. 
Hence  $\hMU^{ev}(*)\cong MU^*$, and for a compact manifold $B$ the group $\hMU(B)$ is an  $MU^*$-module.
We set $\hat R(B):=\hMU^*(B)\otimes_{MU^*}R$ and define the structure maps $R,I,a$ by tensoring the corresponding structure
maps for $\hMU$. Here we identify
$R^*(B)\cong MU^*(B)\otimes_{MU^*}R$ and $\Omega(B,R\R)\cong \Omega(B,MU\R)\otimes_{MU^*}R$. 
The only non-trivial point to show is that the sequence
$$R(B)\stackrel{c_R}{\to} \Omega(B,R\R)/\im(d)\stackrel{a}{\to} \hat R(B)\stackrel{I}{\to} R(B)\to 0$$
is exact. Let us reformulate this as the exactness of 
\begin{equation}\label{eeeqqq}
0\to \Omega(B,R\R)/c_R(R(B))\to  \hat R(B)\to R(B)\to 0\ .
\end{equation}
We start from the exact sequence
$$0\to \Omega(B,MU\R)/c(MU^*(B))\to  \hMU(B)\to MU^*(B)\to 0\ .$$
Tensoring by $R$ gives
\begin{multline*}\Tor_1^{MU^*}(MU^*(B),R)\to
  (\Omega(B,MU\R)/c(MU^*(B)))\otimes_{MU^*}R\\
\to  \hMU(B)\otimes_{MU^*}R\to MU^*(B)\otimes_{MU^*}R\to 0\ .\end{multline*}
Since the tensor product is right exact  we have
$$ (\Omega(B,MU\R)/c(MU^*(B)))\otimes_{MU^*}R\cong  \Omega(B,R\R)/c_R(R(B))\ .$$
We conclude the exactness of (\ref{eeeqqq}) from Landweber's Theorem \ref{landweber} which states that
$\Tor_1^{MU^*}(MU^*(B),R)\cong 0$.
\end{proof}

\subsubsection{}

Let $p\colon V\to A$ be a proper submersion which is smoothly $MU$-oriented (see \ref{smoop}) by $o_p$.
Recall that $\hat R(V)=\hMU(V)\otimes_{MU_*}R$.
\begin{ddd}
We define the push-forward map $p_!\colon \hat R(V)\to \hat R(A)$ by 
$p_!(x\otimes z) :=p_!(x)\otimes z$. 
\end{ddd}
We must show that the push-forward  is well defined.
Let $u\in MU(*)\cong \hMU^{ev}(*)$. We must show that
$p_!(x\cup u)\otimes z=p_!x\otimes uz$.
This indeed follows from the special case of the projection formula Lemma \ref{proje},
$p_!(x\cup u)=p_!(x)\cup u$.
  
{The smooth $MU$-orientation $o_p$ of the proper submersion $p$ gives rise to a form
$A(o_p)\in \Omega(V,R\R$) which we describe in detail in Definition \ref{udidwqdwqdwqdqdwd}.}
The next  theorem states that  the natural and expected  properties of a push-forward hold true. 
\begin{theorem}\label{multi}
The following diagram commutes:
$$\xymatrix{\Omega(V,R\R)/\im(d)\ar[d]^{\int_{V/A}A(o_p)\wedge\dots}\ar[r]^{\mbox{\hspace{0.8cm}}a}&\hat R(V)\ar[d]^{p_!}\ar[r]^I\ar@/^0.5cm/[rr]^R&R^*(V)\ar[d]^{p_!}&\Omega(V,R\R)\ar[d]^{\int_{V/A}A(o_p)\wedge\dots}\\\Omega(A,R\R)/\im(d)\ar[r]_{\mbox{\hspace{0.8cm}}a}&\hat R(A)\ar[r]_I\ar@/_0.5cm/[rr]_R&R^*(A)&\Omega(A,R\R)}
$$
Furthermore, we have the projection formula
$$p_!(p^*x\cup y)=x\cup p_! y\ ,\quad x\in \hat R(A)\ ,\quad y\in \hat R(V)\ .$$
The push-forward is compatible with pull-backs, i.e. for a Cartesian diagram
$$\xymatrix{W\ar[d]^q\ar[r]^F&V\ar[d]^p\\B\ar[r]^f&A}$$ we have
$$q_!\circ F^*=f^*\circ p_!\colon \hat R(V)\to \hat R(B)\ ,$$
where $q$ is smoothly $MU$-oriented by $f^*o_p$. 

If $\colon C\to V$ is a second proper submersion with smooth $MU$-orientation
$o_r$, then the composition $s:=p\circ r$ has the composed orientation $o_s:=o_p\circ o_r$
(see \ref{com2000}), and we have
$$s_!=p_!\circ r_!\colon \hat R(C)\to \hat R(A)\ .$$ 
\end{theorem}
\begin{proof}
This follows immediately by tensoring with $\id_R$ the corresponding results of the push-forward
for the extension of $(MU,c)$. These are  all proven in Section \ref{ghgtg}.
\end{proof}

\begin{kor}
Let $(R_1,g_1)$ and $(R_2,g_2)$ be two Landweber exact formal group laws
  with corresponding cohomology theories $R_i(B):=MU(B)\otimes_{MU^*}R_i$. Let
  $\phi\colon R_1\to R_2$ be a natural transformation of $MU$-modules. Then
  $\phi$ lifts to a natural transformation of smooth cohomology theories as in
  \cite[Definition 1.5]{bs2009} or \cite[Definition
  1.3]{bunke-2007}, $\hat\phi(B):=\id_{\hMU(B)}\otimes \phi$. 

In particular, we have a (multiplicative) smooth complex orientation $\hMU(B)\to \hat K(B)$
    from smooth complex cobordism to smooth K-theory.
\end{kor}

Here, we use again that $\hat K(B)$ is uniquely determined as a multiplicative
extension of $K$-theory \cite{bs2009}.

\section{Normal $G$-structures and  cobordism theories}\label{gzdqw}

\subsection{Representatives of the  stable normal bundle}

\subsubsection{}

In the present paper we construct  geometric models of smooth extensions of cobordism cohomology theories
associated to the families $G(n)$ of classical groups like $U(n)$, $SO(n)$, $Spin(n)$, or $Spin^c(n)$.
We use the notation $MG(B)$ and are in particular interested in the case where $B$ is a smooth manifold.
A cycle for $MG^n(B)$ is a proper smooth map $W\to B$ with a normal $G$-structure such that $\dim(B)-\dim(W)=n$. The relations are given by bordisms. 

Cycles for the smooth extension will have in addition a geometric normal $G$-structure. In order to make a precise definition
we introduce a rather concrete version of the notion of the stable normal bundle.

\subsubsection{}

Let $X$ be a space or manifold. {For $k\in \nat$
   we denote  by $\underline{\R^k}_X$ the (total space of the)} trivial real vector bundle
$X\times \R^k\to X$.  
Let $f\colon A\rightarrow B$ be a smooth map between manifolds.
\begin{ddd}
A representative of the stable normal bundle of $f$ is a real vector bundle
$N\to  {A}$ together with
an exact sequence
$$0\rightarrow TA\xrightarrow{(df,\alpha)} f^*TB\oplus
\underline{\R^k}_A\rightarrow N\rightarrow 0\ ,$$
where we fix only the homotopy class of the projection to $N$.
\end{ddd}

There is a natural notion of an isomorphism of representatives
  of stable normal bundles.
For an integer $l\in\nat$ it is evident how to define the $l$-fold
stabilization $N(l):=N\oplus  \underline{\R^l}_A$ as representative of the
stable normal bundle with corresponding short exact sequence.

\subsubsection{}\label{pullback}
Let $q\colon C\rightarrow B$ be a smooth map which is transversal to $f$.
Then we have a Cartesian diagram
\begin{equation*}
  \begin{CD}
    C\times_BA @>Q>> A\\
    @VVFV @VVfV\\
    C @>q>> B
\end{CD}
\end{equation*}
of manifolds.
If $$0\to  TA\xrightarrow{(df,\alpha)} f^*TB\oplus \underline{\R^k}_A\xrightarrow{u} N\to 0$$ represents the stable
normal bundle of $f$, then we define the pull-back representative of the stable normal
bundle of $F$ by 
\begin{equation*}
0\to T(C\times_BA)\xrightarrow{(dF,\beta)} F^*TC\oplus
\underline{\R^k}_{C\times_BA} \xrightarrow{\gamma} Q^*N \to 0,
\end{equation*}
with $\beta:=Q^*\alpha\circ dQ$ and $\gamma:=Q^*u\circ (F^*dq\oplus
\id_{\underline{\R^k}_{C\times_BA}})$.
Note that $Q^*(N(l))\cong (Q^*N)(l)$.

 \subsubsection{}\label{composition}
We now discuss the stable normal bundle of a composition.
Let $g\colon B\rightarrow C$ a smooth map and \begin{equation}\label{tezwtewzewe}
0\to TB\stackrel{\stackrel{s}{\leftarrow}}{\xrightarrow{(dg,\beta)}} g^*TC\oplus
\underline{\R^l}_B\xrightarrow{v} M\to 0
\end{equation}
be a representative of the stable normal bundle of $g$.
Then we define 
\begin{equation*}
  0\to TA\xrightarrow{(d(gf),\gamma)} (gf)^*TC\oplus
  \underline{\R^l}_A\oplus\underline{\R^k}_A\xrightarrow{w} N\oplus f^*M\to 0
\end{equation*}
as the associated representative of the stable normal bundle of $g\circ
f$. Here $\gamma:=(f^*\beta\circ df,\alpha)$ and $w:=(u\circ  ( f^*s\oplus \id_{\R^k} ) , f^*v\circ
  \pr_{(gf)^*TC\oplus0\oplus \underline{\R^k}})$, where $s$ is the split indicated in (\ref{tezwtewzewe}).
  This split is unique  up to homotopy (since the space of such splits is convex)  
    so that the homotopy class of  $w$ 
is well defined.

 \subsection{$G$-structures and connections on the stable normal bundle}

\subsubsection{}

Let $G$ be a Lie group with a homomorphism $G\to GL(n,\R)$ and consider an $n$-dimensional real vector bundle
$\xi\to X$. 
\begin{ddd}
A $G$-structure on $\xi$ is a pair $(P,\phi)$ of a $G$-principal bundle $P\to X$ and an isomorphism
of vector bundles $\phi\colon P\times_G\R^n\stackrel{\sim}{\to} \xi$.
\end{ddd}

\begin{ddd}
A geometric $G$-structure on $\xi$ is a triple $(P,\phi,\nabla)$, where
$(P,\phi)$ is a $G$-structure and $\nabla$ is a connection on $P$.
\end{ddd}

Note that the trivial bundle $\underline{ \R^n}_X$ has a canonical $G$-structure with $P=X\times G\to X$.

\subsubsection{}\label{familyg}

In order to define a cobordism theory we consider a sequence of groups
$G(n)$,  $n\in \MM$ for an infinite submonoid $\MM$ of
  $(\nat_{\ge 0},+)$  which fit into a chain of commutative diagrams
$$\xymatrix{G(n)\ar[d]\ar[r]&GL(n,\R)\ar[d]\\G(n+k)\ar[r]&GL(n+k,\R)}\ .$$
Typically, $\MM={\nat}$ or $\MM=2\nat$.
This is in particular used in order to define stabilization.
In order to define the multiplicative structure we require in addition
$$\xymatrix{G(n)\times G(m)\ar[d]\ar[r]&GL(n,\R)\times GL(m,\R)\ar[d]\\G(n+m)\ar[r]&GL(n+m,\R)}\ .$$ 
Examples are $O(n)$, $SO(n)$, $Spin(n)$, or $Spin^c(n)$. In the present paper we are in particular interested
in the complex cobordism theory $MU$. In this case we have $\MM=2\nat$ and
we  set
$G(2n)=U(n)$. 
By abuse of notation we will use the symbol $G$ to denote such a family of groups, and by $MG$ the corresponding cobordism theory.

\subsubsection{}
Let $f\colon A\to B$ be a smooth map between manifold.
\begin{ddd}\label{normalg}
 A representative of a normal $G$-structure on $f$ is given by a pair
$(N,P,\phi)$, where $N$ is a representative of the stable normal bundle, and 
$(P,\phi)$ is a $G(n)$-structure on $N$, where $n:=\dim(N)$, $n\in\MM$.\end{ddd}
 For notational
convenience, we write $N$ instead of the short exact sequence with quotient
$N$ which is also contained in the data of a representative of the stable
normal bundle.

\begin{ddd}\label{geomnorg}
 A representative of a geometric normal $G$-structure on $f$ is given by a quadruple
$(N,P,\phi,\nabla)$, where $N$ is a representative of the stable normal bundle
of $f$, and $(P,\phi,\nabla)$ is a geometric $G(n)$-structure on  $N$, where
$n:=\dim(N)$, $n\in\MM$.
\end{ddd}

There are natural notions of isomorphisms of representatives of normal $G$-structures or geometric normal $G$-structures. 
In the following we discuss the operations "stabilization", "pull-back", and "composition" on the level of representatives of normal $G$-structure and geometric normal $G$-structures.
\subsubsection{}

Let $(N,P,\phi)$ be a representative of a normal $G$-structure on $f\colon A\to B$
and consider  $l\in  \MM $. The stabilization $N(l)$ is
$N\oplus \underline{\R^l}_A$. It has 
a canonical $G(n)\times G(l)$-structure with underlying principal bundle  $P\times G(l)\to A$.
We get a $G(n+l)$-structure with the underlying principal bundle $$P(l):=(P\times G(l))\times_{G(n)\times G(l)}G(l+n)\ .$$

\begin{ddd}\label{stabilg}
We define the stabilization  of $(N,P,\phi)$ by
$(N,P,\phi)(l):=(N(l),P(l),\phi(l))$.
\end{ddd}

Let $(N,P,\phi,\nabla)$ is a representative of a geometric normal $G$-structure, then the connection
$\nabla$ induces a connection $\nabla(l)$ on $P(l)$.

\begin{ddd}\label{stabilgg}
We define the stabilization  of $(N,P,\phi,\nabla)$ by
$$(N,P,\phi,\nabla)(l):=(N(l),P(l),\phi(l),\nabla(l))\ .$$
\end{ddd}

\subsubsection{}

We now consider the pull-back and use the notation introduced in \ref{pullback}.
If $(P,\phi)$ is a $G(n)$-structure on
$N$, then $(Q^*P,Q^*\phi)$ is a $G(n)$-structure on $Q^*N$.

\begin{ddd}\label{pullng}
We define the pull-back of a normal $G$-structure by
$$q^*(N,P,\phi):=(Q^*N,Q^*P,Q^*\phi)\ .$$
\end{ddd}

\begin{ddd}\label{pullngg}
We define the pull-back of a geometric normal $G$-structure by
$$q^*(N,P,\phi,\nabla):=(Q^*N,Q^*P,Q^*\phi,Q^*\nabla)\ .$$
\end{ddd}

\subsubsection{}\label{splitg}
We now discuss the composition. Continuing with  the notation of \ref{composition} we consider
$$A\stackrel{f}{\to} B\stackrel{g}{\to} C$$
and  representatives of normal
$G$-structures $(N,P,\phi)$ and $(M,Q,\psi)$ on $f$ and $g$. The sum $N\oplus f^*M$ has a 
natural $G(n)\times G(m)$-structure  with underlying
$G(n)\times G(m)$-bundle $P\times_A f^*Q$, and therefore a $G(n+m)$-structure with underlying
 bundle $$R:=(P\times_Af^*Q) \times_{G(n)\times G(m)}G(n+m)$$ with isomorphism
$\rho\colon R\times_{GL(n+m)}\R^{n+m}\cong N\oplus f^*M$.

\begin{ddd}\label{compog}
We define the composition of representatives of normal $G$-structures by 
$$(M,Q,\psi)\circ (N,P,\phi):=(N\oplus f^*M, R,\rho)\ .$$
\end{ddd}
  
If $\nabla^P$ and $\nabla^Q$ are connections on $P$ and $Q$, then we get an induced connection $\nabla^R$ on $R$.
\begin{ddd}\label{compogg}
We define the composition of representatives of geometric normal $G$-structures by 
$$(M,Q,\psi,\nabla)\circ  (N,P,\phi,\nabla):=(N\oplus f^*M, R,\rho,\nabla^R)\ .$$
\end{ddd}
 
\subsubsection{}
 The following assertions are obvious.
\begin{lem}\label{functg}
\begin{enumerate}
\item 
On the level of representatives of normal $G$-structures or geometric normal $G$-structures,
pull-back and composition commute with stabilization.
\item 
On the level of representatives of normal $G$-structures or geometric normal $G$-structures,
pull-back and composition are functorial.
\item On the level of representatives of normal $G$-structures or geometric normal $G$-structures, pull-back and composition commute with each other. \end{enumerate}
\end{lem}

\subsection{A cycle model for $MG$}

\subsubsection{}
Let us fix a family of groups $G$ and $\MM$ as in \ref{familyg}. 
It determines a multiplicative cohomology theory which is represented by a Thom spectrum $\bMG$.
The map $G(n)\to GL(n,\R)$ induces a map of classifying spaces
$BG(n)\to BGL(n,\R)$. Let $\xi_n\to BG(n)$ denote the pull-back of the universal $\R^n$-bundle. Then for $n\in\MM$
we define  $\bMG^n:=BG(n)^{\xi_n}$, where for a vector bundle $\xi\to X$ we write $X^\xi$ for its  Thom space.   The family of spaces $\bMG^n$, $n\ge 0$,  fits into a spectrum
with structure maps
$$\Sigma^{ d} \bMG^n\cong BG(n)^{\xi_n\oplus
  \underline{\R^{ d }}_{BG(n)}}\to
BG(n+ {d})^{\xi_{n+1}}\cong \bMG^{n+ {d}}\
{,\qquad n,n+d\in\MM}$$
where we use the canonical Cartesian diagram
$$
\begin{CD}
  \xi_{n}\oplus
    \underline{\R^{ {d}}}_{BG(n)} @>>>\xi_{n+{d}}\\ @VVV @VVV \\BG(n) @>>> BG(n+1)\   .
\end{CD}
$$
The ring structure is induced by
\begin{multline*}
  \bMG^n\wedge \bMG^m\cong BG(n)^{\xi_n}\wedge BG(m)^{\xi_m}\cong (BG(n)\times
  BG(m))^{\xi_n\boxplus\xi_m}\\
\to BG(n+m)^{\xi_{n+m}} \cong \bMG^{n+m}\ ,
\end{multline*}

 {for $n,m\in\MM$,} using the canonical Cartesian diagram
$$
\begin{CD}
  \xi_{n}\boxplus \xi_m @>>> \xi_{n+m}\\ @VVV @VVV \\BG(n)\times
    BG(m) @>>>BG(n+m)\ .
\end{CD}
$$
{For $l\notin \MM$ we set $MG^l:=\Sigma^{l-d}MG^d$, where
    $d\le l$ is maximal with $d\in\MM$. The corresponding structure maps and
    multiplication maps are given as suspensions of the maps described above.}

If $A$ is a manifold (or more generally a finite CW-complex), then   the
homotopy theoretic definition of the cobordism cohomology group is
$$\bMG^n(A):=\lim_k [\Sigma^k A_+,\bMG^{n+k}]\ ,$$
where  {the} limit is taken over the stabilization maps
$$[\Sigma^k A_+,\bMG^{n+k}]\to [\Sigma\Sigma^{k} A_+,\Sigma\bMG^{n+k}]\to [\Sigma^{k+1} A_+,\bMG^{n+k+1}]\ ,$$
and $A_+$ is the union of $A$ and a disjoint base point.
Temporarily we use the bold-face notation of the homotopy theoretic definition of the cobordism cohomology theory.
For details we refer to \cite{MR1886843} or \cite{MR0248858}.

\subsubsection{}

We now present a cycle model of the $G$-cobordism theory.
Let $A$ be a smooth manifold.
\begin{ddd}\label{fftr}
A precycle $(p,\nu)$ of degree $n\in\Z$ over $A$ consists of
a smooth map $p\colon W\rightarrow A$ {from a smooth manifold $W$ of dimension $ \dim(W)=\dim(A)-n$}, and a representative $\nu$  of a normal $G$-structure on $p$ (see \ref{normalg}).
A cycle of degree $n\in\Z$ over $A$ is a precycle
$(p,\nu)$ of degree $n$, where $p$ is proper.
\end{ddd}
There is a natural notion of an isomorphism of precycles.
\subsubsection{}

Let $c:=(p,\nu)$ be a precycle over $A$ and $q\colon B\rightarrow A$
be transverse to $p$.
\begin{ddd}
We define the pull-back
$q^*c:=(q^*p,q^*\nu)$, a precycle over $B$.
\end{ddd}
The pull-back is functorial by Lemma \ref{functg}.

\subsubsection{}
We now consider precycles
 $c=(p,\nu)$ over $A$ and $d=(q,\mu)$ over
 $C$ with underlying maps $p\colon B\to A$ and $q\colon A\rightarrow C$. \begin{ddd}\label{compp}
We define the composition
$$d\circ c:=(q\circ p, \mu\circ \nu)$$ using \ref{compog}.
\end{ddd}
The composition $d\circ c$ is a precycle over $C$.
The composition is associative and compatible with pull-back.

\subsubsection{}
Let $c:=(p,\nu)$, $p\colon W\rightarrow A$, and $d:=(q,\mu)$,
$q\colon V\rightarrow B$ be precycles over $A$ and $B$.
Then we can form the diagram
$$
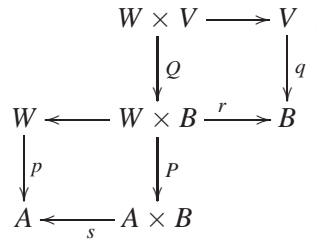
\xymatrix{&W\times V\ar[d]^Q\ar[r]&V\ar[d]^q\\
W\ar[d]^p&W\times B\ar[l]\ar[d]^P\ar[r]^r&B\\
A&A\times B\ar[l]^s&} .
$$


\begin{ddd}\label{prod12}
We define the product of the precycles $c$ and $d$ to be the precycle
$c\times d:=s^*c \circ r^*d$
  over $A\times B$. 
\end{ddd}
Note that
there is an equivalent definition based on
on the diagram
 $$\xymatrix{&W\times V\ar[d]\ar[r]&W\ar[d]^p\\
V\ar[d]^q&A\times V\ar[l]\ar[d]\ar[r]&A\\
B&A\times B\ar[l]&} \ .
$$
It follows from the functoriality of the composition and  its compatibility with the pull-back
that the product of precycles is associative. 

\subsubsection{}
We consider a precycle $b:=((f,p),\nu)$ over $\R\times A$.
\begin{ddd}
The precycle $b$ is called a bordism datum if $f$ is transverse to $\{0\}\in \R$
and $p_{|\{f\ge 0\}}$ is proper. We define the precycle
$\partial b:=i^*b$, where $i\colon A\rightarrow \R\times A$, $i(a):=(0,a)$.
\end{ddd}

\subsubsection{}
Let $c=(p,\nu)$ be a precycle and $l\in\nat$.
\begin{ddd}
We define the $l$-fold stabilization of the precycle $c$
by $c(l):=(p,\nu(l))$ (see \ref{stabilg}).
\end{ddd}

 
 

\subsubsection{}\label{semico}
We now come to the geometric picture of the cobordism theory $MG$.
We consider a smooth manifold $A$ and let $ZMG(A)$ denote the semigroup of isomorphism classes of cycles over $A$ with
respect to disjoint union.  Recall that a relation $\sim$ on a semigroup is compatible with the semigroup structure if
$a\sim b$ implies that $a+c\sim b+c$ for all $c$.
\begin{ddd}\label{geommodel}
Let   ``$\sim$'' be the  minimal equivalence relation
which is compatible with the semigroup structure
and satisfies:
\begin{enumerate}
\item If $b$ is a bordism datum, then $\partial b\sim 0$.
\item If $l\in\nat$, then $c(l)\sim c$.
\end{enumerate} 
We let $MG(A):=ZMG(A)/\sim$ denote the quotient semigroup.
\end{ddd}

\subsubsection{}
Let $0$ denote the cycle of degree $n$  {given by} the empty manifold. 
The following Lemma will be useful in calculations.
\begin{lem}\label{reduce}
Let $c$ be a cycle which is equivalent to $0$.
Then there exists a bordism datum $b$ and $l\ge 0$ such that
$c(l)\cong \partial b$.
\end{lem}
We leave the proof to the interested reader.

\subsubsection{}\label{top1}
We now describe the functoriality, the product, orientations, and the
integration on the level of cycles.
\begin{enumerate}
\item{\bf Functoriality}
Let $f\colon B\to A$ be a smooth map and $x\in MG(A)$. We can represent
$x$ by a cycle $c=(p,\nu)$ such that $p$ and $f$ are transverse.
Then $f^*x$ is represented by $f^*c$.
\item{\bf Product}
Let $c$ and $d$ be cycles for $x\in MG^n(A)$ and $y\in
MG^m(B)$.
Then $x\times y\in MG^{n+m}(A\times B)$ is represented by the cycle
$c\times d$ (see \ref{prod12}). We get the interior product using the pull-back along the diagonal.
\item {\bf Integration} Let $d$ be a cycle over $A$ with underlying map $q\colon V\rightarrow A$.
In this situation we have an integration map $$q_!\colon M {G}(V)\to M {G}(A)\ .$$
 If $x\in MG(V)$ is represented by the cycle
$c$, then $q_!(x)$ is represented by
the cycle $d\circ c$ (see \ref{compp}).
\item{\bf Suspension}
Let $i\colon *\to S^1$ denote the embedding of a point.
For each $d\in\MM$, the trivial bundle $\underline{\R^{ {d}}}_{*}$ represents the stable normal bundle which of course has a canonical 
$G( {d})$-structure.
In this way $i$ is the underlying map of a cycle $\{i\}\in ZMG^1(S^1)$ which represents a class $[i]\in MG^1(S^1)$.

For a manifold $A$ we define $MG^k(A)\rightarrow MG^{k+1}(S^1\times A)$,
 $x\mapsto \{i\}\times x$, which on the level of cycles is represented by 
$c\mapsto \{i\}\times c$. This transformation is essentially the suspension
morphism (not an isomorphism, since we  {neither use reduced
  cohomology nor the suspension of $A$}).
\end{enumerate}

\subsubsection{}
In order to show that the operations defined above  {on the
  cycle level} descend through the equivalence relation $\sim$ the following
observations are useful. Let $b=((f,p),\mu)$ be a bordism datum over $A$ with
underlying map $(f,p)\colon W\to \R\times A$.
Assume that $q\colon B\to A$ is transverse to $p$ and $p_{|\{f=0\}}$.
Then we can form the bordism datum  $(\id_\R\times q)^*b$ over $B$
which will be denoted by $q^*b$.
Note that $$q^*\partial b\cong \partial q^*b\ .$$

Let $e$ be a cycle over $B$. Then we can form $b\times e$
which we can interpret as a bordism datum over $B\times A$.
Note that $$\partial (b\times e)\cong\partial b\times e\ .$$

Let $d$ be a cycle  with underlying map $A\to B$. Let $\pr\colon \R\times B\to B$ be the projection.
Then we can form the bordism datum $\pr^*d\circ b$ over $B$. Note that
$$\partial( \pr^*d\circ b)\cong d\circ \partial b\ .$$

Finally, if   $c$ is a cycle over $W$, then we can form the bordism datum
$b\circ c$ over $B$, and have
$$\partial(b\circ c)\cong \partial b\circ c\ .$$

\subsubsection{}\label{except2}

We now have a geometric and a homotopy theoretic picture of the $G$-cobordism theory
which we distinguish at the moment by using roman and bold-face letters.

\begin{prop} \label{rr5}
There is an isomorphism of ring-valued functors
$MG(A)\cong \bMG(A)$ on compact manifolds.
This isomorphism preserves the product and is compatible with push-forward.
\end{prop}\newcommand{\bMU}{{\mathbf{MU}}} 
\begin{proof}
This follows from the  Pontryagin-Thom construction. 
Since this construction for cobordism cohomology (as opposed to
  homology) seems  {not} to be so well known  let us shortly indicate the main ideas. 
For concreteness let us consider the case of complex cobordism $MU$ and even $ {2}n$.
We have 
$$\bMU^{ {2}n}(A)\cong \colim_i[\Sigma^{ {2}i}A,\bMU^{ {2}n+ {2}i}]\ .$$
Let $h\colon \Sigma^{ {2}i} A\to \bMU^{ {2}n+ {2}i}$ represent some class in $\bMU^{ {2}n}(A)$. Recall that  $\bMU^{ {2}n+ {2}i}=BU(n+ {i})^{\xi_{n+i}}$ is the Thom space the universal bundle $\xi_{n+i}\to BU(n+i)$.
The latter is itself the colimit of Thom spaces $$ BU(n+i)^{\xi_{n+i}}\cong \colim_k  \Gr_{n+i}(\C^{n+i+k})^{\xi_{n+i }} $$ of tautological bundles $\xi_{n+i }$ over the Grassmannians $\Gr_{n+i}(\C^{n+i+k})$ of $\ {(}n+i {)}$-dimensional subspaces in $\C^{n+i+k}$. We can assume that
$h$ factors over some Thom space  
$\Gr_{n+i}(\C^{n+i+k})^{\xi_{n+i}}$, and that the induced map
$S^{ {2}i}\times A\stackrel{p}{\to} \Sigma^{ {2}i}
A\stackrel{f}{\to} \Gr_{n+i}(\C^{n+i+k})^{\xi_{n+i}}$ is smooth  and
transverse to the zero section of $\xi_{n+i}$, where $p$ is the canonical
projection. The preimage of the zero section is a submanifold
$W\subset S^{ {2}i}\times A$ of codimension
$ {2}n+2 i$. We let
$f\colon W\to A$ be induced by the projection.
We use the standard embedding $S^{ {2}i}\to \R^{ {2}i+ {2}}$  in order to trivialize
 the bundle $TS^{ {2}i}\oplus
 \underline{\R}_{S^{ {2}i}}\cong S^{ {2}i}\times
 \R^{ {2}i+ {2}}$.
The embedding $W\hookrightarrow S^{ {2i}}\times A$ thus induces naturally an embedding
$$TW\to T(A\times S^{ i})_{|W}\cong f^*TA\oplus TS^{ {2}i}_{|W}\to \ {f^*TA\oplus}\underline{\R^{ {2}i+ {2}}}_{W}\ .$$
Moreover, the differential of $h$ identifies
the normal bundle
$N:=f^*TA\oplus \underline{\R^{ {2}i+ {2}}}_{W}/TW$
with the pull-back $h_{|W}^*\xi_{n+i} {\oplus\underline\C_W}$, which has a canonical complex structure.
In this way we get the 
normal bundle sequence \begin{equation*}
  0\to T {W}\to f^*TA\oplus\underline{R^{ {2i+2}}}_M\to N\to 0
\end{equation*} and the $U$-structure $\nu=(N,P,\phi)$ 
on $N$.
Note that $f\colon W\to A$ is proper so that we get a cycle
$(f,\nu)$ of degree $n$. One now proceeds as in the case of bordism homology
and shows that the  class $[f,\nu]\in MU^{ {2n}}(A)$ only
depends on the class 
$[h]\in \bMU^{ {2n}}(A)$. 
In this way we get a map $\bMU^{ {2n}}(A)\to MU^{ {2n}}(A)$.

Conversely
one starts with a cycle $(f,\nu)$ of degree $n$.
One observes that up to
stabilization and  homotopy the normal bundle sequence 
\begin{equation*}
  0\to TW\to f^*TA\oplus\underline{R^{ {k}}}_W\to N\to 0
\end{equation*}
  comes from an embedding of $i\colon W\hookrightarrow S^{ {k-1}}\times A$ such that $f=\pr_A\circ i$.
Then we let $ W\to BU(n+ { {(n+k)/2}})$ be a classifying map of $N$
 {(necessarily, $\dim_{\R}N=n+k$ is even)}. It gives rise to a map of Thom spaces
$ W^N\to BU(n+i)^{\xi_{ {(n+k)/2}}}$. We finally precompose with
the clutching map $\Sigma^{ {k}} A\to W^N$ in order to get a map 
$h\colon \Sigma^{ {k}} A\to \bMU^{n+ {k}}$.

One checks that this construction gives the inverse  map $\bMU^n(A)\to
MU^n(A)$. A further standard argument checks that these maps are compatible
with the abelian group and  ring structures and the push-forward. 
\end{proof}

In view of Proposition \ref{rr5} we can drop the bold-face notation for the homotopy theoretic cobordism.
%
%
%
%
%
%
%

\subsection{Power series and genera}

\subsubsection{}

The basic datum for a multiplicative smooth extension of a generalized cohomology theory $h$ is a pair $(h,c)$, where $c\colon h\to HR$ is a natural transformation from $h$ into the ordinary cohomology with coefficients in a graded ring $R$ over $\R$. The transformation $c$ induces
in particular a homomorphism of coefficients $h^*\to R^*$. Our construction of smooth extensions of cobordism theories is based on a description of
$c$ in terms of  characteristic numbers of stable normal bundles.

A ring  homomorphisms $c\colon MG^{{*}}\to R^{{*}}$ is called a $G$-genus. One can classify $SO$ and $U$-genera in terms
of formal power series (see \cite{MR1189136} and \ref{mucase}). Genera for other cobordism theories can be derived from transformations like
$MSpin\to MSO$.  Since the details in the real and complex case differ
slightly, in the present paper we restrict to our main example  {$G:=MU$,
i.e.~$\MM=2\nat_{\ge 0}$, $G(2n)=U(n)$}.
It is easy to modify the constructions for other cases like $MSpin^c$, $MSO$ or $Spin^c$.

\subsubsection{}\label{phi111}

Let $R$ be a commutative $\Z$-graded algebra over $\R$ with $1\in R^0$. 
By $R[[z]]$ we denote the graded ring of formal power series, where $z$ has degree $2$.
Let $\phi\in R[[z]]^0$ be a power series of the form
$1+\phi_1 z +\phi_2 z^2+ \dots$ (note that $\deg(\phi_i)=-2i$). To such a power series we associate
a genus $r_\phi\colon MU^*\rightarrow R^*$  {as in \cite[Section 19]{MS}}.

\begin{theorem}[\cite{MR1189136}] \label{mucase}
The correspondence $\phi\to r_\phi$ gives a bijection between the set $R[[z]]^0$
and $R$-valued $U$-genera.
\end{theorem}
In the following we describe the associated  natural transformation
$r_\phi\colon MU(A)\rightarrow H(A,R)$ of cohomology theories on the level of
cycles, following the procedure as described in \cite{MS}.

\subsubsection{}
We define the power series $K_\phi\in R[[\sigma_1,\sigma_2,\dots,]]^0$ (where
$\sigma_i$ has degree $2i$) such that 
$$K_\phi(\sigma_1,\sigma_2,\dots)=\prod_{i=1}^\infty  \phi(z_i)$$ holds
if we replace $\sigma_i$ by the elementary symmetric functions
$\sigma_i(z_1,\dots)$.

\subsubsection{}
Note that $G(2k)=U(k)$ (see \ref{familyg}).
Let $N\rightarrow W$ be an $n$-dimensional real vector bundle  {for $n$ even} with a $G(n)$-structure
$(P,\phi)$.
Then we have Chern classes $c_j(N):=c_j(P)\in H^{2j}(W,\R)$.
\begin{ddd}
We define the characteristic class $$\phi(N):=K_\phi(c_1(N),c_2(N),\dots)\in H^0(A,R)\ .$$
\end{ddd}

The following properties are well-known (see \cite{MR1189136}).
\begin{lem}\label{propr}
\begin{enumerate}
\item
Let $\underline{\R^{k}}_A$ have the trivial $G(k)$-structure.
Then we  have $\phi(\underline{\R^{k}})_A=1$ for all $k\ge 0$.
\item
If $M$ is a second bundle with a $G(m)$-structure, and $N\oplus M$ has the induced $G(n+m)$-structure, then 
we have $\phi(N\oplus M)=\phi(N)\cup \phi(M)$.
\item If $f\colon B\rightarrow A$ is a continuous map, then we have
$f^*\phi(N)=\phi(f^*N)$, if we equip $f^*N$ with the induced $G(n)$-structure.
\end{enumerate}
\end{lem}

\subsubsection{}

Consider a cycle  $c=(p,\nu)\in ZMU(A)$ of degree $n$ with underlying map $p\colon W\rightarrow A$ and normal $U$-structure $\nu=(N,P,\phi)$.
Then $p$ is a proper map which is oriented for the ordinary cohomology theory $HR$. In particular, we have an integration
$p_!\colon H^*(W,R)\to H^{*{+}n}(A,R)$.
 \begin{ddd}
We define
$$\tilde r_\phi(c):=p_!(\phi(N))\in H^{n}(A,R)\ .$$
\end{ddd}

\subsubsection{}
The following Lemma implies half of Theorem \ref{mucase}.
What is missing is the argument that every $R$-valued $U$-genus comes from
a formal power series.
\begin{lem}\label{mucase1}
The map $\tilde r_\phi$ descends through $\sim$ and induces a natural transformation
$r_\phi\colon MU(A)\to H(A,R)$ of ring-valued functors.
\end{lem} 
\begin{proof}
Using the first and second property in \ref{propr} one checks that
$\tilde r_\phi(c)=\tilde r_\phi(c(l))$.

Assume that $b=((f,q),\mu)$ with underlying map $ {(f,q)\colon
  W\rightarrow \R\times} A$ and $\mu=(M,Q,\lambda)$ is a bordism datum. Then
we get the  {Cartesian} diagram
$$\xymatrix{V\ar[r]^i\ar[d]^p&W\ar[d]^{(f,q)}\\\{0\}\times A\ar[r]&\R\times A}$$ 
 {Set} $N {:=} i^*M$. Therefore
$p_!(\phi(N))= p_!(\phi(i^*M))=p_!(i^* \phi(M))=0$ by the bordism invariance of the push-forward in ordinary cohomology and the third property of \ref{propr}.
Thus the transformation $r_\phi$ is well defined.

It is natural since for $f\colon B\rightarrow A$ which is transverse to $p$
we have a Cartesian diagram
$$\xymatrix{F^*N\ar[d]\ar[r]&N\ar[d]\\f^*V\ar[d]^q\ar[r]^F&V\ar[d]^p\\
B\ar[r]^f&A}\ ,$$ the bundle $F^*N$ represents the stable normal bundle of $q$, and 
$$q_!(\phi( {F}^*N))= q_!(F^*\phi(N))=f^* p_!(\phi(N))$$
by the projection formula.
This implies that
$f^*\tilde r_\phi(c)=\tilde r_\phi(f^*c)$.

 We claim that the transformation is also multiplicative.  
To this end we consider a cycle
$d=(q,\mu)$ with underlying map $q\colon V\rightarrow B$ and normal $G$-structure
$\mu=(M,Q,\lambda)$.
Then the underlying proper map of $c\times d\in ZMU(A \times B)$  
is $p\times q\colon W\times V\rightarrow A\times B$, and the  bundle
$N\boxplus M$ represents its normal $G$-structure.
We thus have
$$(p\times q)_!(\phi(N\boxplus M))=
(p\times q)_!(\phi(N)\times \phi(M))=p_!(\phi(N))\times
q_!(\phi(M))\ .$$
This implies
$$\tilde r_\phi(c\times d)=\tilde r_\phi(c)\times \tilde r_\phi(d)\ .$$
\end{proof}

\subsubsection{}\label{univcase}

The most important example for the present paper is
given by the ring $MU\R:=MU^*\otimes_\Z \R$.
The $MU^*$-module $MU\R$ is Landweber exact. Hence, 
for a compact manifold or finite $CW$-complex $A$ we have
$H^*(A,MU \R)\cong MU^*(A)\otimes_{MU^*}MU\R$
and therefore a  canonical natural transformation
$r\colon MU^*(A)\to H^*(A,MU \R)$, $x\mapsto x\otimes 1$.
This transformation is a genus $r=r_\phi$ for a certain power series
$\phi\in MU \R[[x]]^0$. We refer to \cite{MR1189136} for further details on $\phi$.

\section{The smooth extension of $MU$}\label{wuzwqud}

\subsection{Characteristic forms}

\subsubsection{}

Let $\phi\in R[[z]]^0$ be as in \ref{phi111} and $G$ be the family of groups \ref{familyg} associated to $U(n)$, $n\ge 0$..
We first lift the construction of the characteristic class $\phi(N)\in H^0(A,R)$ of vector bundles $N\to A$
with $G(n)$-structure
 to the form level.

Let $(P,\psi,\nabla^N)$ be a geometric $G(n)$-structure on $N\to A$. 
By $R^{\nabla^N}\in \Omega^2(A,\End(N))$ we denote the curvature of the connection $\nabla^N$.
The fiber-wise polynomial bundle morphism $\det\colon \End(N)\to \underline{\R}_A$ extends 
to $\det\colon  \Omega^{ev}(A,\End(N))\to \Omega^{ev}(A)$.
As usual we define the Chern forms
$c_i(\nabla^N)\in \Omega^{2i}(A)$ by
$$1+c_1(\nabla^N)+c_2(\nabla^N)+\dots=\det(1+\frac{1}{2\pi i}R^{\nabla^N})\ .$$ 
\begin{ddd}
If $N\to A$ is a real vector bundle with a geometric $G(n)$-structure, then we
define $$\phi(\nabla^N):=K_\phi(c_1(\nabla^N),c_2(\nabla^N),\dots)\in \Omega^0(A,R)\ .$$
\end{ddd}

\subsubsection{}

The properties  stated in Lemma \ref{propr} lift to the form
level by well-known properties of the Chern-Weil calculus. 
\begin{lem}\label{propr1}
\begin{enumerate}
\item
Let $k\ge 0$ and $\underline{\R^{k}}_A$ have the trivial $G(k)$-structure with the trivial connection.
Then we  have $\phi(\nabla^{\underline{\R^{k}}_A})=1$.
\item
If $M\to A$ is a second bundle with a geometric $G(m)$-structure and assume that $N\oplus M$ has the induced geometric $G(n+m)$-structure, then 
we have $\phi(\nabla^{N\oplus M})=\phi(\nabla^N) \wedge \phi(\nabla^M)$.
\item  {Assume that} $f\colon B\rightarrow A$ is a smooth map. Then we have
$f^*\phi(\nabla^{N})=\phi(\nabla^{f^*N})$, if we equip $f^*N$ with the induced geometric $G(n)$-structure.
\end{enumerate}
\end{lem}

\subsubsection{}

\begin{ddd}\label{geomcyc}
A geometric precycle over $A$ is a pair $(p,\nu)$ of a smooth map $p\colon V\to A$ and a geometric normal $G$-structure $\nu$ (see \ref{geomnorg}).
A geometric precycle is a cycle if $p$ is proper.
\end{ddd}
Usually we will denote geometric precycles by $\tilde c$, where $c$ denotes the underlying precycle.
Since a principal bundle always admits connections, every precycle can be refined to a geometric precycle.
If $\nu=(N,P,\phi,\nabla)$, then we will write $\nabla^\nu:=\nabla$.

\subsubsection{}

Let $\Omega_{-\infty}(A):=C^{-\infty}(A,\Lambda^*T^*A)$ denote the
differential forms with distributional coefficients. We identify this space
with the topological dual of $C^\infty_c(A,\Lambda^{n-*}T^*A\otimes \Lambda_A)$,
where $\Lambda_A\to A$ is the real orientation  bundle and $n=\dim(A)$. For
this 
identification we use cup product and integration of an $n$-form with values
in the orientation bundle over $A$.
Finally, we define $\Omega_{-\infty}(A,R):=\Omega_{-\infty}(A)\otimes_\R R$ using the algebraic tensor product.

A morphism of complexes inducing an isomorphism in cohomology is called a quasi-isomorphism.
It is well-known (see \cite{MR760450}, or do this exercise using Lemma \ref{choice1}) that the inclusion $\Omega(A)\hookrightarrow \Omega_{-\infty}(A)$ is a quasi-isomorphism. Hence,
$\Omega(A,R)\hookrightarrow \Omega_{-\infty}(A,R)$
is a quasi-isomorphism, too.

\subsubsection{}
Let $p\colon V\to A$ be a proper smooth oriented map. The orientation of  $p$ gives an isomorphism
$p^*\Lambda_A\stackrel{\sim}{\to} \Lambda_V$. We then define the push-forward
$$p_!\colon \Omega_{-\infty}(V)\to \Omega_{-\infty}(A)$$ of degree  $\dim(A)-\dim(V)$
 {by the formula} 
$$<p_!\omega,\sigma>=<\omega,p^*\sigma>\ , \quad \omega\in \Omega_{-\infty}(V)\ ,\sigma\in \Omega(A,\Lambda_A)$$
holds true.
By tensoring with the identity of $R$ we get the map $p_!\colon \Omega_{-\infty}(V,R)\to \Omega_{-\infty}(A,R)$.
 Stokes' theorem implies
$$p_!\circ d=d\circ p_!\ .$$
We get an induced map in cohomology such that the following diagram commutes {:}
\begin{equation}\label{derham}
  \begin{CD}
    H^*(\Omega_{-\infty}(V,R)) @>{deRham}>{\cong}> H^*(V,R)\\
    @VV{p_!}V @VV{p_!}V\\
      H^*(\Omega_{-\infty}(A,R)) @>{deRham}>{\cong}> H^*(A,R)\ .
  \end{CD}
\end{equation}

\subsubsection{}

Let $\tilde c=(p,\nu)$ be a geometric cycle of degree $n$.
\begin{ddd}\label{tcdef}
We define $T(\tilde c):=p_!(\phi(\nabla^\nu))\in \Omega_{-\infty}^{n}(A,R)$.
\end{ddd}
This form is closed, and by (\ref{derham}) we have the following equality in de Rham cohomology:
\begin{equation}\label{voim}
[T(\tilde c) ]=p_!(\phi(N))=\tilde r_\phi(c)\ .
\end{equation}

\subsubsection{}
We now consider a bordism datum $b=((f,q),\mu)$ over a manifold $A$ with 
$(f,q)\colon W\rightarrow \R\times A$.
We build the composition
$$q_!\circ \chi_{\{f\ge 0\}}\colon \Omega^k(W)\rightarrow
\Omega_{-\infty}^{k+l}(A)\ ,$$ where $l:=\dim(A)-\dim(W)$,
and  $\chi_U$ is the multiplication
 operation with the characteristic function of the subset $U$.
Stokes' theorem implies in this case that
\begin{equation}\label{disc1}
d\circ q_!\circ \chi_{\{f\ge 0\}}- 
q_!\circ \chi_{\{f\ge 0\}}\circ d=(q_0)_!\circ i^*\ ,
\end{equation}
where $q_0\colon W_0\to A$ is defined by the Cartesian diagram
$$\begin{array}{ccc}
W_0&\stackrel{i}{\rightarrow}&W\\
q_0\downarrow&&q\downarrow\\
A&\stackrel{a\mapsto (0,a)}{\rightarrow}&\R\times A\end{array}\ ,$$
i.e.~$q_0$ is the underlying map of $\partial b$.

\begin{ddd}
Let $\tilde b:=((f,q),\tilde \nu)$ be a geometric refinement of $b$.
We define $$T(\tilde b):=
q_!\circ \chi_{\{f\ge 0\}}(\phi(\nabla^{\tilde \nu}))\in \Omega_{-\infty}(A) {.}$$
\end{ddd}
Equation (\ref{disc1}) shows that
\begin{equation}\label{bordcurv}d T(\tilde b)=T(\partial \tilde b)\ .\end{equation}

\subsection{The smooth extension of $MU$}  

\subsubsection{}
In the present subsection we construct the smooth extension associated to
the pair $(MU,r_\phi)$, where $\phi\in R[[z]]^0$ is as in \ref{phi111}, and $r_\phi$ is the associated natural transformation
$MU(A)\to H(A,R)$.
Recall the notions of a cycle and a geometric cycle from 
\ref{fftr} and \ref{geomcyc}. The cycles for the smooth extension $\hat{MU}$ of $MU$ will be called smooth cycles.

\begin{ddd}
A smooth cycle of degree $n$ is a pair
$\hat c:=(\tilde c,\alpha)$, where $\tilde c$ is a geometric cycle of degree $n$,
and $\alpha\in \Omega^{n-1}_{-\infty}(A,R)/\im(d)$ is such that 
$$T(\tilde c)-d\alpha:=\Omega(\hat c)\in \Omega^n(A,R)\ .$$
\end{ddd}
The point here is that $T(\tilde c)-d\alpha$ is a smooth representative of
the cohomology class represented by  $T(\tilde c)$. The latter is in general a singular form.
To be explicit note that  in the definition
above $$\im(d):=\im(d\colon \Omega_{-\infty}^{n-2}(A,R)\to
\Omega_{-\infty}^{n-1}(A,R))\ ,$$ i.e.\ we
allow differentials of forms with distribution coefficients.

\subsubsection{}
There is {an evident} notion of an isomorphism of smooth cycles.
We form the graded semigroup $Z\hat{MU}(A)$ of isomorphism classes of  smooth cycles such that the sum
is given by  $$(\tilde c,\alpha)+(\tilde c^\prime,\alpha^\prime)=(\tilde c+\tilde  c^\prime,\alpha+\alpha^\prime)\ ,$$
where, as in the non-geometric case, $\tilde c+\tilde c^\prime$ is given by the disjoint union.
\subsubsection{}

The smooth cobordism group $\hMU(A)$ will be defined as the quotient of $Z\hMU(A)$ by an equivalence relation generated by stabilization and bordism.

\begin{ddd}
Let ``$\sim$'' be the minimal equivalence relation on $Z\hat{MU}(A)$ which is compatible with the semigroup structure (see \ref{semico}) and such
that
\begin{enumerate}
\item For $l\in \MM $ we have
$(\tilde c,\alpha)\sim (\tilde c(l),\alpha)$, where
$\tilde c(l)$ is the $l$-fold stabilization defined by $(p,\nu)(l):=(p,\nu(l))$ (see \ref{stabilgg}).
\item For a geometric bordism datum $\tilde b$ we have
$(\partial \tilde b,T(\tilde b))\sim 0$.
\end{enumerate}
We define  $\hMU^n(A):=Z\hMU^n(A)/\sim$ as the semigroup of equivalence classes of
smooth cycles of degree $n$.
\end{ddd}
We will write $[\tilde c,\alpha]$ for the equivalence class of $(\tilde c,\alpha)$.
 
\subsubsection{}

\begin{lem}
$\hat{MU}^n(A)$ is a group.
\end{lem}
\begin{proof}
Let $[\tilde c,\alpha]\in\hat{MU}(A)$. It suffices to show that it admits an inverse.
Since $MU(A)$ is a group there exists a cycle 
$c^\prime$ such that  
$c + c^\prime \sim  0$. By Lemma \ref{reduce} we can assume that $c(l)+c^\prime(l) \sim \partial b$ for some
bordism datum $b$ and $l\in \nat$. We extend $b$ to a geometric bordism datum $\tilde b$ by choosing a connection
such that $\partial \tilde b\cong \tilde c(l) +\tilde c^\prime(l)$ for some geometric extension $\tilde c^\prime$ of $c^\prime$. Then we have
 $[\tilde c^\prime,T(\tilde b)-\alpha]+[\tilde c,\alpha]=0$. \end{proof}

\subsubsection{}

We now define the structure maps $a,R,I$ (see \ref{ddd556})
of the smooth extension $\hat{MU}$.
\begin{ddd}
\begin{enumerate}
\item
We define $R\colon \hat{MU}(A)\rightarrow \Omega_{d=0}(A,R)$ 
by $R([\tilde c,\alpha]):=T(\tilde c)-d\alpha$.
\item
We define  $a\colon \Omega(A,R)\rightarrow \hMU(A)$ by
$a(\alpha):=[\emptyset,-\alpha]$. 
\item
We define $I\colon \hMU(A)\rightarrow MU(A)$ by
$I([\tilde c,\alpha]):=[c]$ (using the geometric model \ref{geommodel})
\end{enumerate}
\end{ddd}

\begin{lem}
These maps are well defined.
We have $R\circ a=d$.
\end{lem}
\begin{proof}
The only non-obvious part is the fact that $R$ is well defined.
To this end consider a geometric bordism datum $\tilde b$. Then we have
$$R[\partial \tilde b,T(\tilde b)]=T(\partial \tilde b)-dT(\tilde b)=0$$ by
Equation (\ref{bordcurv}). \end{proof}

\subsubsection{}\label{tildepull}
We now extend $A\mapsto \hMU(A)$ to a contra-variant  functor on the category of smooth manifolds.
Let $f\colon B\rightarrow A$ be a smooth map. Then we must construct a
functorial pull-back
$f^*\colon \hMU(A)\rightarrow \hMU(B)$
such that the transformations $R,I,a$ above become natural.

Let $(\tilde c,\alpha)$ be a smooth cycle with $\tilde c=(p,\nu)$, $p\colon W\to A$. We can assume
that $p$ is transverse to $f$.
Otherwise we replace $p$ by a bordant (homotopic) map and correct
$\alpha$ correspondingly {so that the new pair represents the same class in $\hMU(A)$ as $(\tilde c,\alpha)$.}
Then we have the Cartesian diagram
$$\xymatrix{
B\times_AW\ar[d]^P\ar[r]^F&W\ar[d]^p\\
 B\ar[r]^f&A}\ .$$
 The map $P$ is the underlying map of a geometric
 cycle $f^*\tilde c=(P,f^*\nu)$, where $f^*\nu$ is the pull-back of the geometric normal $G$-structure as defined in \ref{pullngg}.
We want to define $f^*[\tilde c,\alpha]:=[f^*\tilde c,f^*\alpha]$.
The problem is that $\alpha$ is a distribution. In order to define the pull-back $f^*\alpha$ of a distributional form 
we need the additional assumption that
$WF(\alpha)\cap N(f)=\emptyset$,
where $N(f)\subseteq T^*A\setminus 0_A$ is the normal set to $f$ given by
$$N(f):=\clo \{\eta\in T^*A\setminus 0_A\:|\:\exists b\in B \:\:s.t.\:\:
f(b)=\pi(\eta)\:\: \mbox{and} \:\:df(b)^*\eta=0\}\ $$
(where $\pi\colon T^*A\rightarrow A$ is the projection),
and $WF(\alpha)$ denotes the wave front set of $\alpha$.  {The
  wave front set of  a distributional form $\alpha$  on $A$ is a conical
  subset of $T^*A$ which measures the locus and the directions of the
  singularities of $\alpha$. For a precise definition and for the properties
  of distributions using the wave front set needed below we refer to
  \cite[ {Section 8}]{MR1996773}.}. 
{Note that we can change $\alpha$ by  exact forms with distribution coefficients without altering the class of $(\tilde c,\alpha)$.}
The idea is to show that one can choose $\alpha$ such that $WF(\alpha)\cap N(f)=\emptyset$ holds.
 {By \cite[Theorem 8.2.4]{MR1996773}, in} this case $f^*\alpha$
is defined.  {It is} independent of the choice again up to exact forms with
distribution coefficients.
The details will be explained in the following paragraphs.

\subsubsection{}

\begin{lem}\label{choice1}
Let $\alpha\in \Omega^n_{-\infty}(A)$.
Then there exists $\beta\in \Omega^{n-1}_{-\infty}(A)$
such that $WF(\alpha-d\beta)\subseteq WF(d\alpha)$.
\end{lem}
\begin{proof}
 We choose a Riemannian metric on $A$. Then we can define the formal adjoint $\delta:=d^*$ of the de Rham differential and the
Laplacian $\Delta:=\delta d+d\delta$. Since $\Delta$ is elliptic we can choose a proper
pseudo-differential parametrix $P$ of $\Delta$.   {This is a pseudo-differential operator of degree $-2$ which is an inverse of $\Delta$ up to pseudo-differential operators of degree $-\infty$ (smoothing operators). A pseudo-differential operator on $A$  is called proper if the restriction of the two projections from the support  (a subset of $A\times A$) of its distribution kernel  to the two
factors $A$ are proper maps.}

Then we form
$G:=\delta P$. This pseudo-differential operator satisfies
$d G+G d = 1 + S$, where $S$ is a proper smoothing operator.
We thus can set $\beta := G\alpha$ and have
$$\alpha-d\beta = G d\alpha - S\alpha \ .$$
Since $S\alpha$ is smooth and $WF(Gd\alpha)\subseteq WF(d\alpha)$ (a pseudo-differential
operator does not increase wave front sets)
we see that
$WF(\alpha-d\beta)\subseteq WF(d\alpha)$. \end{proof} 
If $\alpha\in \Omega_{-\infty}(A,R)$, then we can write for some $s\in \nat$
$$\alpha=\sum_{i=1}^s \alpha_i\otimes r_i$$ with $\alpha_i\in \Omega^n_{-\infty}$,
and with linearly independent $r_i\in R$. In this case the wave front set of $\alpha$ is by definition
$WF(\alpha):=\cup_{i=1}^s WF(\alpha_i)$. It is now easy to see that Lemma \ref{choice1} extends to forms with coefficients in $R$.

\subsubsection{}
 
\begin{lem}\label{choice}
If $(\tilde c,\alpha)$, $\tilde c=(p,\nu)$, is a smooth cycle, then we can
choose $\alpha$ such that $WF(\alpha)\subseteq N(p)$.
\end{lem}
\begin{proof}
It is a general fact that the wave front set of the push-forward of a smooth distribution along a map is contained in the normal set of the map. In view of  {Definition} \ref{tcdef} we have
$WF(T(\tilde c))\subseteq N(p)$.
Since $T(\tilde c)-d\alpha$ is smooth we have
$WF(d\alpha)=WF(T(\tilde c))\subseteq N(p)$, and by Lemma \ref{choice1} we can change $\alpha$
by an exact form such that $WF(\alpha)\subseteq N(p)$. \end{proof}

\subsubsection{}

A reformulation of the fact that $f$ and $p$ are transverse is that
$N(f)\cap N(p)=\emptyset$. Using Lemma \ref{choice} we now take
a representative of $\alpha$ such that $WF(\alpha)\subseteq N(p)$.
Then $f^*\alpha$ is a well defined distribution.
\begin{ddd}
We define $f^*[\tilde c,\alpha]=[f^*\tilde c,f^*\alpha]$,
where we take representatives $\tilde c=(p,\nu)$  and $\alpha$
such that $f$ and $p$ are transverse and $WF(\alpha)\subseteq N(p)$.
\end{ddd}

\subsubsection{}
\begin{lem}\label{pulli}
The pull-back is well defined and functorial.
\end{lem}
\begin{proof}
First we show that the pull-back is well defined with respect to the choice of $\alpha$.
Let $\beta\in \Omega_{-\infty}(A,R)$ and $\alpha^\prime:=\alpha+\beta$ be such that
$T(\tilde c)-d\alpha^\prime$ is smooth. This implies that
$WF(\alpha^\prime)\subseteq N(p)$, and hence 
 $WF(d\beta)\subseteq N(p)$. By Lemma \ref{choice1} we can  modify $\beta$ by a closed form  
 such that $WF(\beta)\subseteq N(p)$.
Then $f^*\alpha^\prime=f^*\alpha+d f^*\beta$.  

It is easy to see that the pull-back is additive and preserves stabilization.
It remains to show that it preserves zero bordism.
Let $\tilde b=((h,q),\mu)$ be a geometric bordism datum over $A$ with $(h,q)\colon W\to \R\times A$. We define $W_0:=h^{-1}(\{0\})$ and assume  that
$q$ and $q_{|W_0}$ are  transverse to $f$. 
We then have the geometric bordism datum $(\id_\R\times f)^*\tilde b$ over $B$.

Let us define the normal datum of $b$ by
\begin{eqnarray*}N(b)&:=&\clo\{\eta\in T^*A\setminus 0_A\:|\: \exists v\in W
\:\:s.t.\:\: E(v)=\pi(\eta)\:\:
\mbox{and} \:\:\left[dE(v)^*\eta=0\right.\\&&\left.\mbox{or}\:\: v\in W_0 \:\: \mbox{and} \:\:dE(v)^*\eta_{|T_vW_0}=0\right]\} .\end{eqnarray*}
Then we have $WF(T(\tilde b))\subseteq N(b)$.
Again, since $q$ and $q_{|W_0}$ are  transverse to $f$ we have $N(b)\cap
N(f)=\emptyset$ so that $f^*T(\tilde b)$ is well defined.
Using the fact that in a Cartesian diagram push-forward of distributions commutes with pull-back  we get $f^*T(\tilde b)=T(f^*\tilde b)$. 
It follows that
$(f^*\partial \tilde b ,f^*T(\tilde b))=(\partial f^*\tilde b,T(f^*\tilde b))$.
This implies that the pull-back is well defined on the level of equivalence classes.

We now show functoriality. Let $g\colon C\to B$ be a second smooth map.
If $\hat x\in  \hMU(A)$, then we can choose the representing smooth cycle
$(\tilde c,\alpha)$ with  $\tilde c=(p,\nu)$ such that $p$ is transverse to
$f$ and $f\circ g$. In this case one  easily sees that
$(f\circ g)^*(\tilde c, \alpha)$ and $g^*f^*(\tilde c,\alpha)$ are
isomorphic cycles.
\end{proof}

\subsubsection{}
We now have defined a functor $A\mapsto \hMU(A)$ from smooth manifolds to graded groups.
\begin{lem}
The transformations $R$, $I$ and $a$ are natural.
\end{lem}
\begin{proof}
Straightforward.\end{proof}

\subsubsection{}

We now define the outer product
$$\times \colon \hMU(A)\otimes \hMU(B)\rightarrow \hMU(A\times
B)\ .$$

Let $\hat x\in \hMU(A)$ be represented by $(\tilde c,\alpha)$, and let
$\hat y\in \hMU(B)$ be represented by $(\tilde e,\beta)$.
In \ref{prod12} we have already defined the product of cycles $c\times
e$. Here we enhance this definition to the geometric level.
Write $\tilde c=(p,\nu)$ and $\tilde d=(q,\mu)$. Then we define
$\tilde c\times \tilde d:=(p\times q,\nu\oplus \mu)$, where the sum of geometric normal $G$ structures $\nu\oplus \mu$
is defined similarly as in the non-geometric case.

Note that we have a graded outer product
$$\times\colon \Omega_{-\infty}(A,R)\otimes \Omega_{-\infty}(B,R)\rightarrow
\Omega_{-\infty}(A\times B,R)\ .$$
\begin{ddd}
We define the product of smooth cycles $(\tilde c,\alpha)\times (\tilde e,\beta)$ by 
$$(\tilde c\times \tilde e,(-1)^{|\hat x|}R(\hat x)\times \beta + \alpha\times
T(e))\ ,$$ and we define
the product $\hat x\times  \hat y\in  \hMU(A\times B)$ to be the corresponding equivalence class.
 \end{ddd}
This cycle level definition needs a few verifications.

\begin{lem}\label{outpro}
\begin{enumerate}
\item
The outer product is well defined.
\item
It is associative, i.e. $(\hat x\times\hat y)\times \hat z=\hat
x\times (\hat y\times \hat z)$, where $\hat z\in \hMU(C)$.
\item It is graded commutative in the sense that
$F^*(\hat x\times \hat y)=(-1)^{|\hat x||\hat y|}\hat y\times \hat x$,
where $F\colon B\times A\rightarrow A\times B$ is the flip $F(b,a):=(a,b)$.
\item The product is natural, i.e. if
$f\colon C\to A$ is a smooth map, then we have $f^*\hat x\times \hat y=(f\times
\id_B)^*(\hat x\times \hat y)$.
\end{enumerate}
\end{lem}
\begin{proof}
We first show that the cycle level definition of the outer product passes through the equivalence relation.
It is obvious that the outer product is bilinear and preserves
stabilizations in both arguments.  It remains to verify that it preserves zero bordisms.

Let $\tilde b$ be a geometric bordism datum. 
Then we can form the geometric bordism datum $\tilde b\times \tilde e$ (see \ref{prod12}).
We have $T(\tilde b\times \tilde e)=T(\tilde b)\times T(\tilde e)$ so that
\begin{eqnarray*}(\partial \tilde b,T(\tilde b))\times (\tilde e,\beta)&=&(\partial \tilde b\times \tilde e,T(\tilde b)\times T(\tilde e))\\&=&(\partial (\tilde b\times\tilde
e),T(\tilde b\times \tilde e))\\&\sim& 0\ .\end{eqnarray*}

In order to see that the product also preserves zero bordism in the
second entry we rewrite
\begin{equation}\label{eweu}
(-1)^{|\hat x|}R(\hat x)\times \beta + \alpha\times
T(\tilde e)\stackrel{\im(d)}{\equiv} (-1)^{|\hat x|} T(\tilde c)\times \beta +
\alpha\times R(\hat y)
\end{equation}
and apply the same argument as above.
Associativity, graded commutativity, and naturality hold true
on the level of smooth cycles. {To see this, for} commutativity we use again (\ref{eweu}), and the proof of 
associativity is  based on similar calculations. 
\end{proof}

\subsubsection{}
As usual, the outer product determines a graded commutative ring structure by restriction to the diagonal.
\begin{ddd}
We define the ring structure on $\hMU(A)$ by
$\hat x\cup \hat y:=\Delta^*(\hat x\times \hat y)$,
where $\Delta\colon A\rightarrow A\times A$ is the diagonal.
\end{ddd}

The following assertions are consequences of Lemmas \ref{pulli} and  \ref{outpro}.
\begin{kor}
$A\mapsto \hMU(A)$ is a contra-variant functor from the category of  manifolds to 
the category of graded  {commutative} rings.

\end{kor}

\begin{lem}
The transformations $R$ and $I$ are multiplicative, and we have $a(\alpha)\cup \hat x=a(\alpha \wedge R(\hat x))$.
\end{lem}
\begin{proof}
Straightforward calculation. \end{proof}

\subsubsection{}

Recall that we have fixed {in \ref{mucase}, \ref{mucase1},
  \ref{univcase}} a graded ring $R$ over $\R$ and a formal power series
$\phi\in R[[z]]^0$ which determines an $R$-valued  $U$-genus $r_\phi$.
\begin{theorem}\label{fundprop}
The functor $\hMU$ together with the transformations
$R,I,a$ is a multiplicative smooth extension of the pair $(MU,r_\phi)$.
\end{theorem}
\begin{proof}
We must verify the properties required in Definitions \ref{ddd556} and \ref{multdef1}.
Most of them have been shown above. We are left with the commutativity of
\begin{equation}\label{ghgh67}
\xymatrix{\hMU(B)\ar[d]^R\ar[r]^I& MU(B)\ar[d]^{r_\phi}\\
\Omega_{d=0}(B,R)\ar[r]^{dR}&H(B,R)}\ .
\end{equation}
and the exactness of 
\begin{equation}\label{exax1}
MU(B)\stackrel{r_\phi}{\to} \Omega(B,R)/\im(d)\stackrel{a}{\to} \hMU(B)\stackrel{I}{\to} MU(B)\to 0\ .
\end{equation}
The commutativity of the diagram (\ref{ghgh67}) is a direct consequence of (\ref{voim}).

We now discuss exactness of (\ref{exax1}).
We start with the surjectivity of $I$. Let $x\in MU(B)$ be represented by a cycle $c$.
Then we can choose a geometric refinement $\tilde c$.
We have $dT(\tilde c)=0$, and by Lemma \ref{choice1} there
exists $\alpha\in \Omega_{-\infty}(B,R)$ such that $T(\tilde c)-d\alpha$ is
smooth. Therefore $(\tilde c,\alpha)$ is a smooth cycle, and we have $x=I[\tilde c,\alpha]$.

We now discuss exactness at $\hMU(B)$.
It is clear that $I\circ a=0$.
Let $\hat x\in  \hMU(B)$, be such that $I(\hat x)=0$.
Then we can assume that $\hat x$ is of the form $[\partial \tilde b,\alpha]$
for some geometric bordism datum $\tilde b$. Hence
$\hat x=a(T(\tilde b)-\alpha)$.

We now show exactness at $\Omega(B,R)/\im(d)$.
Let $x\in MU(B)$ be represented by a cycle $c$.
Then we choose a geometric refinement $\tilde c$, and by \ref{choice1} a form
$\alpha\in \Omega_{-\infty}(B,R)$ such that $T(\tilde c)-d\alpha$ is smooth.
We have $r_\phi(x)=T(\tilde c)-d\alpha$. Let $c=(p,\nu)$ with $p\colon V\to B$, and consider the constant map $h\colon V\to \R$ with value $1$. The geometric normal $U$-structure of $(h,p)\colon V\to \R\times B$ can also  be represented by  $\nu$.
Then $\tilde b=((h,p),\nu)$ is a geometric bordism datum with $\partial \tilde b=\emptyset$ and $T(\tilde b)=T(\tilde c)$.
It follows
$$a(d\alpha-T(\tilde c))=[\partial \tilde b,T(\tilde b)-d\alpha]=[\partial \tilde b,T(\tilde b)]=0\ .$$
This proves that $a\circ r_\phi=0$.

Let now $\alpha\in \Omega(A,R)$ be such that
$a(\alpha)=0$.
Then there exist geometric bordism data $\tilde b_0, \tilde b_1$ such that
$\partial \tilde b_0\cong \partial \tilde b_1$ and 
$T(\tilde b_0)-T(\tilde b_1)- \alpha\in\im(d)$.
This already implies that $\alpha$ is closed.
We construct a geometric cycle $\tilde c$ 
such that $T(\tilde c)=T(\tilde b_0)-T(\tilde b_1)$ by glueing the bordism
data along their common boundary.  
Then $[\alpha]=[T(\tilde c)]=r_\phi([c])$ in de Rham cohomology.
\end{proof}

\subsection{Smooth $MU$-orientations}

\subsubsection{}

As before we fix a graded ring $R$ over $\R$ and a formal power series
$\phi\in R[[z]]^0$. Let $\hMU$ be the smooth extension of $(MU,r_\phi)$ as in Theorem \ref{fundprop} with structure maps $R,a,I$.
If $q\colon V\to A$ is a proper $MU$-oriented map, then we have an integration
$q_!\colon MU(V)\to MU(A)$ (see \ref{top1}).
Under the assumption that $q$ is a submersion we introduce the notion of a smooth $MU$-orientation 
and define the integration map $q_!\colon \hMU(V)\to \hMU(A)$.
\subsubsection{}

Let $q\colon V\to A$ be a proper submersion.
\begin{ddd}
A representative of a smooth $MU$-orientation of $q$ is a pair
$\check{c}:=(\tilde c,\sigma)$, where $\tilde c$ is a geometric cycle with underlying map $q\colon V\to A$ and
$\sigma\in \Omega^{-1}(V,R)$. 
\end{ddd}
A representative of a smooth $MU$-orientation of $q$ induces in particular an $MU$-orientation of $q$.

\subsubsection{}

We now introduce an equivalence relation $\sim$ called stable homotopy on the set of representatives of smooth $MU$-orientations of $q$.
 \begin{ddd}
We define the $l$-fold stabilization of  
$(\tilde c,\sigma)$ by $(\tilde c,\sigma)(l):=(\tilde c(l),\sigma)$.
\end{ddd}

Let $h_i\colon A\to \R\times A$  denote the inclusions $h_i(a):=(i,a)$, $i=0,1$.
Consider a geometric cycle $\tilde d=(p,\mu)$ over $\R\times A$ with underlying map
$p:=\id_{\R}\times q\colon \R\times V\to \R\times A$. It gives rise to a closed form
$\phi(\nabla^\mu)\in \Omega^0(\R\times V,R)$.
Let $\tilde c_i:=h_i^*\tilde d$, $\tilde c_i=(q,\nu_i)$.
\begin{ddd}
We call $\tilde d$ a homotopy between $\tilde c_0$ and $\tilde c_1$.
\end{ddd}

\begin{ddd}
We define the transgression form
$$\tilde \phi(\nabla^{\nu_1},\nabla^{\nu_0}):=\int_{[0,1]\times V/V}\phi(\nabla^\mu)\in \Omega^{-1}(A,R)/\im(d)\ .$$
\end{ddd}
Since the underlying cycle $d$ of $\tilde d$ is a product, and since the space of geometric refinements of $d$ is contractible,
the transgression form is well defined  {independent of} the choice of the homotopy (this is a standard argument in the theory of characteristic forms). By Stokes' theorem the transgression satisfies
\begin{equation}\label{transg}
d\tilde \phi(\nabla^{\nu_1},\nabla^{\nu_0})=\phi(\nabla^{\nu_1})-\phi(\nabla^{\nu_0})\ .
\end{equation}

\begin{ddd}
We call two representatives of a smooth $MU$-orientation $(\tilde c_i,\sigma_i)$
homotopic if there exists a  homotopy $\tilde d$ from $\tilde c_0$ to $\tilde c_1$, and
$\sigma_1-\sigma_0=\tilde \phi(\nabla^{\nu_1},\nabla^{\nu_0})$.
\end{ddd}

\subsubsection{}

We now define equivalence of representatives of smooth $MU$-orientations.
\begin{ddd}\label{smoop}
Let $\sim$ be the minimal equivalence relation on the set of representatives of smooth $MU$-orientations on $q$ such that
\begin{enumerate}
\item $(\tilde c,\sigma)\sim (\tilde c(l),\sigma)$
\item $(\tilde c_0,\sigma_0)\sim (\tilde c_1, \sigma_1)$, if $(\tilde c_0,\sigma_0)$ and $(\tilde c_1, \sigma_1)$ are homotopic.
\end{enumerate}
A smooth $MU$-orientation of  $q$ is an equivalence class of representatives of smooth $MU$-orientations which we will usually write as $o:=[\tilde c,\sigma]$.
\end{ddd} 

\subsubsection{}
Let $\tilde c:=(q,\nu)$ and $\check{c}:=(\tilde c,\sigma)$ be a representative of a smooth $MU$-orientation. 
\begin{ddd}
We define $A(\check{c}):=\phi(\nabla^\nu)-d\sigma\in \Omega^0(V,R)$.
\end{ddd}
\begin{lem}\label{udidwqdwqdwqdqdwd}
The form $A(\check c)$ only depends on the smooth $MU$-orientation $[\check{c}]$ represented by $\check{c}$.
\end{lem}
\begin{proof}
This immediately follows from (\ref{transg}) and the definition of homotopy. \end{proof} 
Below we will write $A(o):=A(\check{c})$, where $o:=[\check{c}]$.

\subsubsection{}

In the following two paragraphs we define the operations of pull-back and composition of smooth $MU$-orientations.
We start with the pull-back.
Let $f\colon B\to A$ be a smooth map which is transverse to $q$. Then we have the Cartesian diagram
$$\xymatrix{W\ar[d]^Q\ar[r]^F&V\ar[d]^q\\
B\ar[r]^f&A}\ .$$
\begin{ddd}\label{bdwn}
We define the pull-back of a representative  of a smooth $MU$-orientation of $q$ by
$f^*(\tilde c,\sigma):=(f^*\tilde c,F^*\sigma)$ (see \ref{tildepull}) which is a representative of a smooth $MU$-orientation of $Q$.
\end{ddd}

\begin{lem}
The pull-back is compatible with the equivalence relation. It induces a functorial
pull-back of smooth $MU$-orientations. We have $A(f^*o)=F^*A(o)$. 
\end{lem}
\begin{proof}
It is clear that the pull-back is compatible with stabilization.
Let $\tilde d$ be a homotopy from $\tilde c_0$ to $\tilde c_1$.
Then $(\id_{\R}\times f)^*\tilde d$ is a homotopy from
$f^*\tilde c_0$ to $f^*\tilde c_1$. Furthermore, one checks that
$\tilde \phi(\nabla^{f^*\nu_1},\nabla^{f^*\nu_0})=f^*\tilde \phi(\nabla^{\nu_1},\nabla^{\nu_0})$.
These formulas imply that the pull-back preserves homotopic representatives of smooth $MU$-orientations.
We conclude that the pull-back is well defined on the level equivalence
classes. Functoriality and the  {fact that $A(f^*o)=F^*A(o)$ are} easy to see.
\end{proof} 

\subsubsection{}

We now define the composition {of smooth $MU$-orientations}.
Let $p\colon A\to B$ be a second proper submersion, and let $(\tilde d,\theta)$, $\tilde d=(p,\mu)$, be a representative of a smooth
$MU$-orientation of $p$. Let $o_q=[\tilde c,\sigma]$ and $o_p:=[\tilde d,\theta]$. By  $\tilde d\circ \tilde c$ we denote the composition of geometric cycles which is based on Definition \ref{compogg}.

\begin{ddd}\label{com2000}
We define 
$$o_p\circ o_q:=[\tilde d\circ \tilde c,A(o_q)\wedge q^*\theta+\sigma\wedge q^*\phi(\nabla^{\mu})]\ .$$
\end{ddd}
The definition requires some verifications.
\begin{lem}
The composition of smooth $MU$-orientations is well defined, compatible with pull-back, and functorial.
\end{lem}
 \begin{proof}
We first show that the composition is well defined.
It is clear that the composition is compatible with stabilization.
Let $\tilde b$ be a homotopy from $\tilde c_0$ to $\tilde c_1$.
Then $\pr_2^*\tilde d\circ  \tilde b$ is a homotopy from $\tilde b\circ\tilde c_0$ to $\tilde b\circ\tilde c_1$, where $\pr_2:\R\times B\to B$ is the projection.
We further calculate (using the properties stated in Lemma \ref{propr1})
\begin{eqnarray*}
(\sigma_1-\sigma_0) \wedge q^*\phi(\nabla^{\mu})&=&
 \tilde \phi(\nabla^{\nu_1},\nabla^{\nu_0}) \wedge q^*\phi(\nabla^{\mu})  \\
&=&\tilde \phi(\nabla^{\mu\circ\nu_1},\nabla^{\mu\circ\nu_0})\ .
\end{eqnarray*}
This calculation implies that the composition $(\tilde d,\theta)\circ\dots $ preserves
homotopic representatives.

Let us now consider a homotopy
$\tilde e$ from $\tilde d_0$ to $\tilde d_1$
We get a homotopy
$\tilde e\circ \tilde c$ from $\tilde d_0\circ\tilde c$ to $\tilde d_1\circ\tilde c$.
Furthermore we rewrite (note that we work modulo $\im(d)$)
$$A(o_q)\wedge q^*\theta+\sigma\wedge q^*\phi(\nabla^{\mu})=\phi(\nabla^\nu)\wedge q^*\theta + \sigma\wedge q^*A(o_p)\ .$$
We have 
\begin{eqnarray*} \phi(\nabla^\nu)\wedge q^*(\theta_1-\theta_0)&= &\phi(\nabla^\nu)\wedge q^*\tilde \phi(\nabla^{\mu_1},\nabla^{\mu_0})\\
&=&\tilde \phi(\nabla^{\mu_1\circ \nu},\nabla^{\mu_0\circ\nu})\ .
\end{eqnarray*}
Hence $\dots \circ (\tilde c,\nu)$ preserves homotopic representatives.
This finishes the proof that the composition is well defined.

\subsubsection{}
The composition of smooth $MU$-orientations is associative and compatible with pull-back. For completeness let us state the second fact in greater detail. Let $r\colon Q\to B$ be a map which is transverse to
$q$ and $p\circ q$. Then we have the composition of pull-back diagrams
$$\xymatrix{Q\times_BV\ar[r]\ar[d]&V\ar[d]^q\\Q\times_BA\ar[d]\ar[r]^s&A\ar[d]^p\\Q\ar[r]^r&B}\ .$$
 {In this situation} we have
$$s^*o_p\circ r^*o_q=r^*(o_p\circ o_q)\ .$$
 {We leave the} details of the straightforward proof to the reader. \end{proof}

\subsection{The push-forward}\label{ghgtg}

\subsubsection{}

Let $p\colon V\to A$ be a proper submersion with a smooth $MU$-orientation $o_p:=[\tilde d,\sigma]$, $\tilde d=(p,\nu)$.
 In the following, $(\tilde c,\alpha)$ denotes a smooth cycle on $V$, and we use the notation
$$\int_{V/A}:=p_!:\Omega_{-\infty}(V,R)\to \Omega_{-\infty}(A,R)$$ for the integration of forms.
\begin{ddd}\label{cpush}
We define the push-forward on the level of cycles by
$$p_!(\tilde c,\alpha)=(\tilde d\circ  \tilde c,\int_{V/A}(\phi(\nabla^{\nu})\wedge \alpha+\sigma\wedge R(\tilde c,\alpha)))\ .$$
 \end{ddd}
\begin{lem}
For fixed $(\tilde d,\sigma)$ the push-forward preserves equivalence of smooth cycles.
Furthermore, the induced map
$p_!\colon \hMU(V)\to \hMU(A)$ only depends on the equivalence class $[\tilde d,\sigma]$ of representatives of the smooth $MU$-orientation.
\end{lem}
\begin{proof}
It is clear that the push-forward is additive and compatible with stabilization. 
Let now $\tilde b$ be a geometric bordism datum over $V$.
Let $\pr\colon \R\times A\to A$ be the projection and form $(\tilde e,\theta):=\pr^*(\tilde d,\sigma)$.
Then $\tilde e\circ \tilde b$ is a bordism datum, and we have
$T(\tilde e\circ \tilde b)=\int_{V/A} \phi(\nabla^{\nu})\wedge T(\tilde b)$.
We calculate
\begin{eqnarray*}
p_!(\partial \tilde b,T(\tilde b))&=&(\tilde d\circ \partial \tilde b,\int_{V/A}\phi(\nabla^{\nu})\wedge T(\tilde b))\\
&=&(\partial (\tilde e\circ \tilde b), T(\tilde e\circ \tilde b))\ .
\end{eqnarray*}
This equality implies that $p_!$ preserves zero bordisms.

For a fixed  representative $(\tilde d,\sigma)$ of the smooth $MU$-orientation  we now have a well defined map
$p_!\colon \hMU(V)\to \hMU(A)$. Next we show that it only depends on the smooth orientation
represented by $(\tilde d,\sigma)$.
Again it is clear that stabilization of the representative of the smooth orientation does not
change $p_!$. We now consider a homotopy $\tilde b$  from $(\tilde d_0,\sigma_0)$ to $(\tilde d_1,\sigma_1)$.
The idea of the argument is to translate this homotopy into a bordism datum.
To this end we first consider a model case.
Let $\kappa\colon \R\to \R$ be defined by $\kappa(x):=x-x^2$. 
Then $\kappa^{-1}(\{[0,\infty)\})=[0,1]$. We choose
a representative of the stable normal bundle of $\kappa$ with a geometric 
$U$-structure $\mu$ such that
$\tilde \kappa=(\kappa,\mu)$ is a geometric bordism datum.

Let $\pr_1\colon \R\times A\to \R$ denote the projection.  The composition $\tilde r:=\pr_1^*\tilde \kappa \circ\tilde b$ is now a bordism datum. 
Let $\rho$ denote the representative of the geometric $U$-structure on the normal bundle
of $r$. We consider $\tilde r\circ  \pr_2^*\tilde c$ as a geometric bordism datum
with $\partial (\tilde r\circ \pr_2^*\tilde c)=\tilde d_0\circ\tilde c +(\tilde d_1\circ\tilde c)^{op}$, where $(\cdot)^{op}$ indicates a flip of the orientation.
\ {Fix} $\tilde c=(q,\nu)$  {with} $q\colon U\to V$ and $\tilde d_i=(p,\lambda_i)$.
\begin{eqnarray*}
T(\tilde r\circ \pr_2^*\tilde c )&=&\int_{q^{-1}r^{-1}([0,\infty)\times V)/A} \phi(\nabla^\rho)\wedge \phi(\nabla^\nu)\\&=&\int_{V/A}\left(\tilde \phi(\nabla^{\lambda_1},\nabla^{\lambda_0})\wedge \int_{U/V}\phi(\nabla^\nu) \right)
\end{eqnarray*}
On the other hand
\begin{eqnarray*}\lefteqn{
\int_{V/A}(\phi(\nabla^{\lambda_1})-\phi(\nabla^{\lambda_0}))\wedge \alpha + (\sigma_1-\sigma_0)\wedge R(\tilde c,\alpha)}&&\\
&=&\int_{V/A}d\tilde \phi (\nabla^{\lambda_1},\nabla^{\lambda_0})\wedge \alpha +  \tilde \phi (\nabla^{\lambda_1},\nabla^{\lambda_0})\wedge R(\tilde c,\alpha)\\
&=&\int_{V/A}\tilde \phi (\nabla^{\lambda_1},\nabla^{\lambda_0})\wedge d\alpha +  \tilde \phi (\nabla^{\lambda_1},\nabla^{\lambda_0})\wedge R(\tilde c,\alpha)\\
&=&\int_{V/A}\left(\tilde \phi (\nabla^{\lambda_1},\nabla^{\lambda_0})\wedge \int_{U/V}\phi(\nabla^{\nu})\right)
\end{eqnarray*}
These two equations together show that
$(\tilde d_1,\sigma_1)\circ(\tilde c,\alpha)\sim (d_0,\sigma_0)\circ (\tilde c,\alpha)$.
Indeed
$$(\tilde d_0\circ\tilde c + (\tilde d_1\circ\tilde c)^{op},\int_{V/A}(\phi(\nabla^{\lambda_1})-\phi(\nabla^{\lambda_0}))\wedge \alpha + (\sigma_1-\sigma_0)\wedge R(\tilde c,\alpha))=
(\partial (\tilde r\circ \pr_2^*\tilde c) ,T(\tilde r\circ\pr_2^*\tilde c))\ .$$
\end{proof} 

\subsubsection{}

\begin{lem}\label{hjhj}
The following diagram commutes.
$$\xymatrix{\Omega(V,R)/\im(d)\ar[d]^{\int_{V/A}A(o_p)\wedge\dots}\ar[r]^{\mbox{\hspace{0.6cm}}a}&\hMU(V)\ar[d]^{p_!}\ar[r]^I\ar@/^0.5cm/[rr]^R&MU(V)\ar[d]^{p_!}&\Omega(V,R)\ar[d]^{\int_{V/A}A(o_p)\wedge\dots}\\\Omega(A,R)/\im(d)\ar[r]_{\mbox{\hspace{0.6cm}}a}&\hMU(A)\ar[r]_I\ar@/_0.5cm/[rr]_R&MU(A)&\Omega(A,R)}
$$\end{lem}
\begin{proof}
Commutativity of the left square
follows from partial integration
$$\int_{V/A} (\phi(\nabla^\nu)\wedge \alpha-\sigma\wedge d\alpha)=\int_{V/A} (\phi(\nabla^\nu)-d\sigma)\wedge  \alpha=\int_{V/A} A(o_p)\wedge\alpha \ .$$
For the right square we use 
$$T(\tilde d\circ \tilde c)=\int_{V/A} \phi(\nabla^\nu)\wedge T(\tilde c) ,$$ which implies
\begin{eqnarray*}R(p_!(\tilde c,\alpha))&=&T(\tilde d\circ \tilde c)-d\int_{V/A}(\phi(\nabla^{\nu})\wedge \alpha+\sigma\wedge R(\tilde c,\alpha))\\
&=&\int_{V/A} (\phi(\nabla^\nu)\wedge T(\tilde c)-\phi(\nabla^{\nu})\wedge d\alpha-d\sigma\wedge R(\tilde c,\alpha) )\\&=&\int_{V/A} (\phi(\nabla^\nu)-d\sigma) \wedge R(\tilde c,\alpha) \\&=&
\int_{V/A} A(o_p)\wedge R(\tilde c,\alpha)\ .
\end{eqnarray*}
Commutativity of the middle square is a direct consequence of geometric description
of $p_!\colon MU(V)\to MU(A)$ (see \ref{top1}).
\end{proof}

\subsubsection{}
Let $p\colon V\to A$ be as before with the smooth $MU$-orientation
$o_p:=[\tilde d,\sigma]$. We furthermore consider a proper submersion $q\colon A\to B$ with a
smooth $MU$-orientation $o_q:=[\tilde e,\rho]$, $\tilde e=(q,\mu)$. Let $r:=q\circ p\colon V\to B$ be equipped with the composed
smooth $MU$-orientation $o_r:=o_q\circ o_p$ (see Definition \ref{com2000})
\begin{lem}
The push-forward is functorial, i.e.\ we have the equality
$$r_!=q_!\circ p_!\colon \hMU(V)\to \hMU(B)\ .$$
\end{lem}
\begin{proof} 
The equality holds on the smooth cycle level.
The proof is a straightforward calculation of both sides by inserting the definitions
and using the right square in Lemma \ref{hjhj}. \end{proof} 

\subsubsection{}

Let $p\colon V\to A$ be a  proper smoothly $MU$-oriented map as above, and let 
$f\colon B\to A$ be a second smooth map so that we get a Cartesian diagram
$$\xymatrix{W\ar[d]^P\ar[r]^F&V\ar[d]^p\\
B\ar[r]^f&A}\ .$$
The map $P$ has an induced smooth $MU$-orientation $o_P:=f^*o_p$ (see Definition \ref{bdwn}).
\begin{lem}\label{carte}
The push-forward commutes with pull-back, i.e. we have the equality
$$P_!\circ F^*=f^*\circ p_!\colon \hMU(V)\to \hMU(B)\ .$$
\end{lem}
\begin{proof}
The equality holds true on the level of smooth cycles
$(\tilde c,\alpha)$ whose underlying map is transverse to $F$.
By definition we have  $o_P=[f^*\tilde d,f^*\sigma]$. Furthermore, it follows immediately
from the definitions that $f^*(\tilde d\circ \tilde c)=f^*\tilde d\circ F^*\tilde c$.
The final ingredient of the verification is the identity
$$f^*\circ \int_{V/A}{\dots} =\int_{W/B}\circ F^*\dots\ .$$  \end{proof}

\subsubsection{}

Let $p\colon V\to A$ be a smoothly $MU$-oriented proper submersion  {as above}.
\begin{lem}\label{proje}
The projection formula holds true, i.e. for $x\in \hMU(A)$ and $y\in \hMU(V)$ we have
$p_!(p^*x\cup y)=x\cup p_!y$.
\end{lem}
\begin{proof}
We consider the diagram
$$\xymatrix{V\ar[d]^p\ar[r]^{(p,\id_)}\ar@/^0.8cm/[rr]^{\id_V}&A\times V\ar[r]^{\pr_2}\ar[d]^q&V\ar[d]^p\\
A\ar[r]^{\Delta}\ar@/_0.8cm/[rr]^{\id_A}&A\times A\ar[r]^{\pr_2}&A}\ ,$$
where $q:=\id_A\times p$ has the induced orientation $o_q:=\pr_2^*o_p$.
If we show that
\begin{equation}\label{eqsss}
q_!(x\times y)=x\times p_!(y)\ ,
\end{equation}
then by the definition of the cup-product and applying Lemma  \ref{carte}  to the left Cartesian square
we get the result. Equation (\ref{eqsss}) holds true on the level of smooth cycles and is straightforward to check. 
\end{proof}



\end{document}